\documentclass{amsart}
\usepackage{amsmath,amssymb, amsbsy}
\usepackage[dvips]{graphicx}
\newcommand{\R}{{\mathbb R}}
\newcommand{\N}{{\mathbb N}}
\newcommand{\Z}{{\mathbb Z}}
\newcommand{\weakly}{\rightharpoonup}
\renewcommand{\a }{\alpha }

\renewcommand{\d }{\delta }
\newcommand{\D }{\Delta }
\newcommand{\e }{\varepsilon }

\newcommand{\n }{\nabla }

\newcommand{\Di}{{\mathcal D}^{1,2}(\R^N)}
\newcommand{\alchi}{\raisebox{1.7pt}{$\chi$}}

\newcommand{\sgn}{\mathop{\rm sgn}}

\newenvironment{pf}{\noindent{\sc Proof}.\enspace}{\hfill\qed
\medskip}
\newenvironment{pfn}[1]{\noindent{\bf Proof of {#1}.\enspace}}{\hfill\qed
\medskip}
\newtheorem{Theorem}{Theorem}[section]
\newtheorem{Corollary}[Theorem]{Corollary}

\newtheorem{Lemma}[Theorem]{Lemma}
\newtheorem{Proposition}[Theorem]{Proposition}

\newtheorem{remark}[Theorem]{Remark}
\begin{document}

\title[Schr\"odinger operators with multipolar inverse-square
potentials]{On Schr\"odinger operators with multipolar
inverse-square potentials} 

\author[Veronica Felli]{Veronica Felli}
\address{\hbox{\parbox{5.7in}{\medskip\noindent{Veronica Felli:
        Universit\`a di Milano Bicocca, Dipartimento di Statistica,
        Via Bicocca degli Arcimboldi 8, 20126 Milano, Italy.
        \em{E-mail address: }{\tt veronica.felli@unimib.it.}}}}}
\author[Elsa M. Marchini \and Susanna Terracini]{Elsa M. Marchini \and
  Susanna Terracini}
\address{\hbox{\parbox{5.7in}{\medskip\noindent{Elsa M. Marchini,
        Susanna Terracini: Universit\`a di Milano Bicocca,
        Dipartimento di Ma\-t\-ema\-ti\-ca e Applicazioni, Via Cozzi 53, 20125
        Milano, Italy. \em{E-mail addresses: }{\tt
          elsa.marchini@unimib.it, susanna.terracini@unimib.it}.}}}}

\date{May 29, 2007}

\thanks{Supported by Italy MIUR, national project
``Variational Methods and Nonlinear Differential Equations''.
}

\begin{abstract}
\noindent Positivity, essential self-adjointness, and spectral properties of a
class of Schr\"odinger 
operators with multipolar inverse-square potentials are discussed. In
particular a necessary and sufficient condition on the masses of
singularities for the existence of at least a configuration of poles
ensuring the positivity of the associated quadratic form is
established. 
\end{abstract}

\maketitle
\tableofcontents
\section{Introduction and statement of the main results}\label{intro}

This paper deals with a class of Schr\"odinger operators associated with potentials possessing multiple inverse square singularities:
\begin{equation}\label{eq:17}
L_{\lambda_1,\dots,\lambda_k,a_1,\dots,a_k}:=-\D
-{\displaystyle{\sum_{i=1}^k}}\dfrac{\lambda_i}{|x-a_i|^2}
\end{equation}
where $N\geq 3$, $k\in\N$, $(\lambda_1,\lambda_2,\dots,\lambda_k)\in\R^k$,
$(a_1,a_2,\dots, a_k)\in \R^{k N}$, $a_i\neq a_j$ for $i\neq j$.

From the mathematical point of view, the main reason of interest in inverse square potentials of type $V(x)\sim{\lambda}{|x|^{-2}}$ relies in their criticality: indeed they have the same homogeneity as the laplacian and do not belong to the Kato's class, hence they cannot be regarded as a lower order perturbation term. Besides, potentials with this rate of decay are critical also in nonrelativistic quantum mechanics, as they represent an intermediate threshold between regular potentials (for which there are ordinary stationary states) and singular potentials (for which the energy is not lower-bounded and the particle falls to the center), for more details see \cite{FLS, LL}.
We also mention that inverse square singular potentials arise in many
other physical contexts: molecular physics, see e.g. \cite{leblond},
quantum cosmology \cite{BE}, linearization of combustion models
\cite{BG, GP, VZ}. Moreover we emphasize the correspondence between
nonrelativistic Schr\"odinger operators with inverse square potentials
and relativistic Schr\"odinger operators with Coulomb potentials, see
\cite{DL}.

In recent literature, several papers are concerned with Schr\"odinger equations with Hardy-type singular potentials, see e.g. \cite{AFP, CCP, CTV, egnell, FG, GP, Jan,  RW,  SM, terracini}.

The case of multi-polar Hardy-type potentials was considered in \cite{duyckaerts, FT, FT2}. More precisely, in \cite{duyckaerts} estimates on resolvant truncated at high frequencies are proved for Schr\"odinger operators with multiple inverse square singular potentials. 
In \cite{FT, FT2} the existence of ground states for a  class of multi-polar nonlinear elliptic equations with critical power-nonlinearity is investigated. 

The purpose of the present paper is to analyze how the mutual interaction among the poles influences the spectral properties of the class of operators \eqref{eq:17} and which configurations of singularities ensure positivity of the associated quadratic form 
$$
Q_{\lambda_1,\dots,\lambda_k,a_1,\dots,a_k}(u):=\int_{\R^N} |\n  u(x)|^2dx-\sum_{i=1}^k
{\lambda_i}\int_{\R^N}\frac{u^2(x)}{|x-a_i|^2}\,dx.
$$
As a natural setting to study the operators defined in (\ref{eq:17}) and the associated quadratic forms, we introduce the functional space $\Di$ defined as the completion of $C^\infty_{\rm c}(\R^N)$ with respect to the Dirichlet norm
$$
\|u\|_{\Di}:=\bigg(\int_{\R^N}|\n u(x)|^2\,dx\bigg)^{1/2}.
$$ 
We introduce now a notion of positivity which is closely related to
the property of {\em{strong subcriticality}} introduced in \cite{DS}, see
also Remark \ref{r:ds} for a discussion about the relations between
these notions in the framework of multi-polar inverse square
potentials. We say that $Q_{\lambda_1,\dots,\lambda_k,a_1,\dots,a_k}$ is {\em positive
  semidefinite} whenever 
$$
Q_{\lambda_1,\dots,\lambda_k,a_1,\dots,a_k}(u)\geq 0,\quad\text{ for all
}u\in\Di, 
$$
whereas it is said to be {\em positive definite} if there exists a positive constant $\e\!=\!\e(\lambda_1,\dots,\lambda_k,a_1,\dots,a_k)$ such that
\begin{equation*}
Q_{\lambda_1,\dots,\lambda_k,a_1,\dots,a_k}(u)\geq \e(\lambda_1,\dots,\lambda_k,a_1,\dots,a_k)\int_{\R^N}|\n u(x)|^2\,dx,\quad\text{ for all }u\in\Di.
\end{equation*}
In the case of a single pole operator $-\D-\frac{\lambda}{|x|^2}$, a complete answer to the question of positivity is provided by the classical Hardy inequality  (see for instance \cite{ GP, HLP}):
\begin{equation}\label{eq:hardy}
\Big(\frac{N-2}2\Big)^2\int_{\R^N}\frac{|u(x)|^2}{|x|^2}\,dx\leq
\int_{\R^N}|\nabla u(x)|^2\,dx,\quad\text{ for all }u\in\Di,
\end{equation}
where the constant $\big(\frac{N-2}{2}\big)^2$ is optimal and not attained. 
The optimality of such a constant implies that the quadratic form $Q_{\lambda,0}$ is positive definite in $\Di$ if and only if $\lambda<\big(\frac{N-2}{2}\big)^2$. 
For a more detailed discussion about the properties of monopole singular Hardy type operators, we refer to \cite{GP,terracini,VZ}. 

As observed in \cite{FT}, the positivity of
$Q_{\lambda_1,\dots,\lambda_k,a_1,\dots,a_k}$ depends on the strength and the
location of the singularities, and more precisely on the shape of the
configuration of poles, due to scaling properties of the operator. In
particular in \cite[Proposition 1.2]{FT} it is proved that a
sufficient condition for the quadratic form to be positive definite
for any choice of $a_1,a_2,\dots,a_k$ is that
$\sum_{i=1}^k\lambda_i^+<\frac{(N-2)^2}4$, where $t^+:=\max\{t,0\}$. 
Conversely, if $\sum_{i=1}^k\lambda_i^+>\frac{(N-2)^2}4$, then it is
possible to find a configuration of poles  such that the quadratic
form is not positive definite. As a consequence, in the case $k=2$, if

$$
\lambda_i<\frac{(N-2)^2}4,\quad\mbox{ for }i=1,2\quad\mbox{ and }\quad\lambda_1+\lambda_2<\frac{(N-2)^2}4,
$$ 
then for any choice of $a_1,a_2$, the quadratic form $Q_{\lambda_1,\lambda_2,a_1,a_2}$ is positive definite. 

The first result of the present paper relies in a necessary and sufficient
condition on the masses $\lambda_i$ to have positivity of the quadratic form $Q_{\lambda_1,\dots,\lambda_k,a_1,\dots,a_k}$ for at least a configuration of poles.

\begin{Theorem}\label{t:mainresult}
Let $(\lambda_1,\dots,\lambda_k)\in\R^k$. Then 
\begin{equation}\label{eq:19}
\lambda_i<\frac{(N-2)^2}4,\quad\mbox{ for every }i=1,\dots,k,\quad\mbox{
and }\quad\sum_{i=1}^k\lambda_i<\frac{(N-2)^2}4 ,
\end{equation}
is a necessary and sufficient condition for the existence of a configuration of poles  $\{a_1,\dots,a_k\}$ such that the quadratic form $Q_{\lambda_1,\dots,\lambda_k,a_1,\dots,a_k}$ associated to the operator $L_{\lambda_1,\dots,\lambda_k,a_1,\dots,a_k}$ is positive definite.
\end{Theorem}

The necessity of condition \eqref{eq:19} follows quite directly from the optimality of the best constant in Hardy's inequality and proper interaction estimates (see the proof in Section \ref{mainresult}). 
To prove the sufficiency we study the possibility of obtaining a coercive operator by summing up multisingular potentials which give rise to positive 
quadratic forms, after pushing them very far away from each other to weaken the
interactions among poles. Since the potentials overlap at infinity, the singularity of the resulting potential is the sum of their masses, so that we need to require a control on it.
To this aim, we consider the following class of potentials
\begin{align*}
{\mathcal V}:=\left\{
\begin{array}{ccc}
V(x)={\displaystyle{\sum_{i=1}^k\frac{\lambda_i \alchi_{B(a_i,r_i)}(x)}{|x-a_i|^2}
+\frac{\lambda_{\infty} \alchi_{\R^N\setminus B(0,R)}(x)}{|x|^2}+W(x):\
k\in\N,\ r_i,R\in\R^+,}}\\[20pt] 
a_i\in\R^N,\ a_i\neq a_j\ \text{for }i\neq j,\
\lambda_i,\lambda_{\infty}\in\big(-\infty,\frac{(N-2)^2}4\big),\ W\in L^{N/2}(\R^N)\cap L^{\infty}(\R^N)
\end{array}
\right\}.
\end{align*}
By Hardy's and Sobolev's inequalities, it follows that, for any
$V\in\mathcal V$, the first eigenvalue $\mu(V)$ of the operator
$-\D-V$ in $\Di$ is finite, namely 
$$
\mu(V)=\inf_{u\in\Di\setminus\{0\}}\frac{\displaystyle\int_{\R^N}
\big(|\n  u(x)|^2-V(x)u^2(x)\big)\,dx}
{\displaystyle\int_{\R^N}|\n u(x)|^2\,dx}\,>-\infty. 
$$
Hence, we shall frame in the class $\mathcal V$ the analysis of
coercivity conditions for Schr\"odinger operators.  
Let us notice that $\mu(V)$ can be estimated from above as follows.

\begin{Lemma}\label{l:estmu}
For any $V(x)=\sum_{i=1}^{k}\lambda_i\alchi_{B(a_i,r_i)}(x)|x-a_i|^{-2}+
 \lambda_{\infty}\alchi_{\R^N\setminus
 B(0,R_0)}(x)|x|^{-2}+W(x)
\in{\mathcal V}$, there holds:
\begin{align*}
i)\quad&\text{if $\lambda_i<0$ for all
$i=1,\dots,k,\infty$, then }\mu(V)\leq 1;\\
ii)\quad &\text{if $\max_{i=1,\dots,k,\infty}\lambda_i>0$, then }\mu(V)\leq
  1-\frac4{(N-2)^2}\,\max_{i=1,\dots,k,\infty}\lambda_i.
\end{align*}
\end{Lemma}

A first positivity result in the class
${\mathcal V}$ relies in the following Shattering Lemma yielding positivity in the case of
 singularities  localized strictly near the
poles.

\begin{Lemma}{\bf [\,Shattering of singularities\,]}\label{l:sep}
  For any $\{a_1,a_2,\dots,a_k\}\subset B(0,R_0)\subset\R^N$, $a_i\neq
  a_j$ for $i\neq j$,
  $\{\lambda_1,\dots,\lambda_k,\lambda_{\infty}\}\subset(-\infty,(N-2)^2/4)$, 
and $0<\alpha<1-\frac{4\max_{i=1,\dots,k,\infty}\lambda_i}{(N-2)^2}$,
there
  exists $\delta=\delta(\alpha)>0$ such that 
\begin{align*}
\mu\bigg(\sum_{i=1}^k\frac{\lambda_i\,\alchi_{B(a_i,\delta)}(x)}{|x-a_i|^2}&+
\frac{\lambda_{\infty}\alchi_{\R^N\setminus B(0,R_0)}(x)}{|x|^2}\bigg)\\
&\geq
\begin{cases}
1-\frac{4\max_{i=1,\dots,k,\infty}\lambda_i}{(N-2)^2}-\alpha,&\text{if
}\max_{i=1,\dots,k,\infty}\lambda_i>0,\\
1,&\text{if }\max_{i=1,\dots,k,\infty}\lambda_i\leq0.
\end{cases}
\end{align*}
In particular, the quadratic form associated to the
  operator
$$
{\mathcal L}^{\delta}_{\lambda_1,\dots,\lambda_k,\lambda_\infty,a_1,\dots,a_k}:=-\Delta
-\sum_{i=1}^k\frac{\lambda_i\,\alchi_{B(a_i,\delta)}(x)}{|x-a_i|^2}-
\frac{\lambda_{\infty}\alchi_{\R^N\setminus B(0,R_0)}(x)}{|x|^2}
$$
is positive definite.
\end{Lemma}

The above result can be extended to Schr\"odinger operators whose
potentials have infinitely many 
singularities localized in sufficiently small neighborhoods of 
equidistanced poles, see Lemma~\ref{l:ret}.

Lemma \ref{l:sep} implies that  Schr\"odinger operators with
potentials in ${\mathcal V}$ are compact perturbations of
positive operators, as stated in the following lemma.

\begin{Lemma}\label{l:semibounded}
For any $V\!\in\!\mathcal V$,
there exist $\widetilde V\!\in\!{\mathcal V}$ and $\widetilde W\!\in\!
L^{N/2}(\R^N)\cap L^{\infty}(\R^N)$ such that   $\mu(\widetilde V)>0$ and
$V(x)=\widetilde
V(x)+\widetilde W(x)$. 
\end{Lemma}

The proof of Theorem \ref{t:mainresult} is based on the following
result, which  yields a powerful tool to obtain positive
operators by choosing properly the configuration of poles. 

\begin{Theorem}{\bf [\,Separation Theorem\,]}\label{t:scattering} 
Let 
\begin{align*}
&V_1(x)=\sum_{i=1}^{k_1}\frac{\lambda_i^1\alchi_{B(a_i^1,r_i^1)}(x)}{|x-a_i^1|^2}+
\frac{\lambda_{\infty}^1\alchi_{\R^N\setminus
B(0,R_1)}(x)}{|x|^2}+W_1(x)\in{\mathcal V},\\ 
&V_2(x)=\sum_{i=1}^{k_2}\frac{\lambda_i^2\alchi_{B(a_i^2,r_i^2)}(x)}{|x-a_i^2|^2}+
\frac{\lambda_{\infty}^2\alchi_{\R^N\setminus
B(0,R_2)}(x)}{|x|^2}+W_2(x)\in{\mathcal V}. 
\end{align*}
Assume that $\mu(V_1),\mu(V_2)>0$, namely that the quadratic forms
associated to the operators $-\Delta-V_1$, $-\Delta-V_2$ are positive
definite, and that $\lambda_{\infty}^1+\lambda_{\infty}^2<(\frac{N-2}{2})^2$.   
Then, there exists $R>0$ such that, for every $y\in\R^N$ with $|y|\geq
R$, the quadratic form associated to the operator
$-\Delta-\left(V_1+V_2(\cdot-y)\right)$ is positive definite.  
\end{Theorem}

\begin{remark}\label{r:ds}
We mention that the notion of positivity (respectively nullity) in the
class ${\mathcal V}$ is related to that of {\it subcriticality}
(respectively {\it criticality}) of potentials arising in the
classification given by Simon \cite{Simon}.  A recent breakthrough in
the theory by Pinchover and Tintarev \cite{pinchovertintarev} relates
subcriticality with the presence of a {\it spectral gap}. More
precisely, the subcritical case occurs if and only if the
corresponding Hardy type inequality can be improved, namely if there
exists $p\in C^0(\R^N\setminus\{a_1,\dots,a_k\})$, $p>0$ in
$\R^N\setminus\{a_1,\dots,a_k\}$, such that
\begin{equation}\label{eq:59}
\nu_p(V)=\inf_{u\in C^{\infty}_{\rm
    c}(\R^N\setminus\{a_1,\dots,a_k\})}\frac{\int_{\R^N}|\n
u(x)|^2\,dx-\int_{\R^N}V(x)\,u^2(x)\,dx}
{\int_{\R^N}p(x)\,u^2(x)\,dx}>0.
\end{equation}
As observed in Remark \ref{rem:segno}, being any operator with
potential in ${\mathcal V}$ compact perturbation of a positive
operator, $\nu_p(V)$ has the same sign as $\mu(V)$.  Therefore in the
class ${\mathcal V}$, subcriticality is equivalent to positivity.
\end{remark}

Concerning the Separation Theorem \ref{t:scattering}, we recall that
the case of potentials with compact support was investigated in
\cite{Simon}, while the case of potentials in the Kato class was
studied in \cite{pinchover95}. We notice that in both these papers the
potential is a lower order perturbation of the laplacian and none of
them include the case of potentials with singularities of inverse
square type as Theorem \ref{t:scattering} does.  Such a case is not a
trivial issue and requires some additional assumptions to control the
singularity at infinity, as one can understand just summing up two
mono-polar operators and observing that in this case the positivity of
the resulting operator does not depend on the distance between poles,
due to scaling properties.

Besides the sign, a natural question arising in the study of $\mu(V)$
is its attainability. Indeed, while in the case of a single pole the
best constant in the associated Hardy's inequality 
is not achieved, when dealing with multipolar
Hardy-type potentials a balance between positive and negative
interactions between the poles can lead, in some cases, to
 attainability of the best constant in the corresponding 
Hardy-type inequality, as proved in the proposition below. 
We mention that in the literature analogous phenomena are studied by
\cite{BM, BMS, MMP} for generalized Hardy's inequalities in bounded
and unbounded domains with a boundary, and by \cite{tertikas} for
potentials satisfying Hardy type inequalities and perturbations of
them in $\R^N$.
\begin{Proposition}\label{p:attai}
Let $V(x)=\sum_{i=1}^{k}\lambda_i\alchi_{B(a_i,r_i)}(x)|x-a_i|^{-2}+
 \lambda_{\infty}\alchi_{\R^N\setminus
 B(0,R_0)}(x)|x|^{-2}+W(x)
\in{\mathcal V}$. If 
\begin{equation}\label{eq:77}
\mu(V)<1-\frac4{(N-2)^2}\,\max\big\{0,\lambda_1,\dots,\lambda_k,\lambda_{\infty}\big\}
\end{equation}
then $\mu(V)$ is attained.
\end{Proposition}

The properties of ${\mathcal V}$ proved in Lemmas \ref{l:sep} and
\ref{l:semibounded}, make such a class of potentials to be a quite
natural setting to study the spectral properties of multisingular
Schr\"odinger operators in $L^2(\R^N)$. Indeed, as a direct
consequence of Lemma \ref{l:semibounded}, we obtain that if
$V=\widetilde V+\widetilde W\in{\mathcal V}$, with $\mu(\widetilde
V)>0$ and $\widetilde W\in L^{N/2}(\R^N)\cap L^{\infty}(\R^N)$, then
   
\begin{equation}\label{eq:20}
\nu_1(V):=\inf_{u\in
H^1(\R^N)\setminus\{0\}}\frac{\displaystyle\int_{\R^N}\big(|\n
u(x)|^2-V(x)u^2(x)\big)\,dx}{\displaystyle
\int_{\R^N}|u(x)|^2\,dx}\geq - \|\widetilde
W\|_{L^{\infty}(\R^N)}>-\infty,
\end{equation}
i.e. Schr\"odinger operators with
potentials in ${\mathcal V}$ are semi-bounded.

In order to ensure {\em{quantum completeness}} of the quantum system
associated to the Schr\"odinger operator, a further key aspect to be
discussed is the {\em essential self-adjointness},
namely the existence of a unique self-adjoint extension. Semi-bounded
Schr\"odinger operators are essentially self-adjoint whenever the
potential is not too singular (see \cite{reedsimon}).  
On the other hand, when dealing with inverse square potentials, the
singularity is quite strong and essential self-adjointness is a nontrivial issue. 
In the case of one pole singularity, it was proved in \cite{kwss} (see Theorem
\ref{t:kwss}) that essential self-adjointness  depends on the value
of the mass of the singularity with respect to the threshold
$(N-2)^2/4-1$. The following theorem extends such a result  to potentials lying in
the class ${\mathcal V}$.

\begin{Theorem}\label{t:self-adjo}
Let 
$V(x)=\sum_{i=1}^k\lambda_i|x-a_i|^{-2}\alchi_{B(a_i,r_i)}(x)+
\lambda_{\infty}|x|^{-2}\alchi_{\R^N\setminus
B(0,R)}(x)+W(x)\in{\mathcal V}$. 
Then $-\D-V$  is essentially self-adjoint
in $C^{\infty}_{\rm c}(\R^N\setminus\{a_1,\dots,a_k\})$ if and only if
$\lambda_i\leq (N-2)^2/4-1$ for all $i=1,\dots,k$. 
\end{Theorem}

The proof of the above theorem is based on the regularity results
proved in \cite{FS3, murata, pinchover94} for elliptic equations with
singular weights (see also \cite{DV,Gutierrez}), which allow us to
give the exact asymptotic behavior near the poles of solutions to
Schr\"odinger equations with singular potentials. The characterization
of essential self-adjointness for multi-singular Schr\"odinger
operators given above can be extended to the case of infinitely many
singularities distributed on reticular structures, see Theorem
\ref{t:self-adjo-ret}.

From Theorem \ref{t:self-adjo}, we have that if
$V\in{\mathcal V}$ with $\lambda_i\leq(N-2)^2/4-1$ for all
$i=1,\dots,k$, then the associated Schr\"odinger operator 
 $-\D-V$ is essentially self-adjoint and, consequently, admits
a unique self-adjoint extension, which is given by the {\em {Friedrichs
extension}} $(-\D-V)^F$:
\begin{align}\label{eq:52}
D\big((-\D-V)^F\big)&=
\{u\in H^1(\R^N):\ -\D u-Vu\in L^2(\R^N)\},\quad
u\mapsto-\D u-Vu.
\end{align}

Otherwise, i.e. if $\lambda_i>(N-2)^2/4-1$ for some $i$, $-\D-V$ is
not essentially self-adjoint and admits  many self-adjoint extensions,
among which the Friedrichs extension is the only one whose domain is
included in $H^1(\R^N)$, namely in the domain of the associated quadratic form (see also \cite[Remark 2.5]{duyckaerts}).
A complete description of the spectrum of the Friedrichs extension of operators
with potentials in the class ${\mathcal V}$ is given in the
proposition below.
\begin{Proposition}\label{p:spe}
For any $V\in\mathcal V$, there holds:
\begin{align*}
{\bf \ref{p:spe}.1.} \ &\text{the essential spectrum }\sigma_{\rm
  ess}\big((-\D-V)^F\big)=[0,+\infty);\\
{\bf \ref{p:spe}.2.}\ &\text{if }\nu_1(V)<0 \text{ then  the discrete spectrum of
}(-\D-V)^F\text{ consists in a finite number of}\\
&\text{negative eigenvalues}.
\end{align*}
\end{Proposition}
The nature of the bottom of the essential spectrum of the Friedrichs
extension of operators of type (\ref{eq:17}) is analyzed in the
following theorem. 
\begin{Theorem}\label{t:bottom}
Let $(\lambda_1,\dots,\lambda_k)\in\R^{k}$ satisfy (\ref{eq:19}).
Then 
\begin{equation}\label{eq:72}
\sum_{i=1}^k\lambda_i^+>\frac{(N-2)^2}4\quad\mbox{
and }\quad\sum_{i=1}^k\lambda_i<\frac{(N-2)^2}4-1,
\end{equation}
is a necessary and sufficient condition for the existence of at least a
configuration of singularities  $(a_1,\dots,a_k)\in\R^{Nk}$,
$a_i\neq a_j$ for $i\neq j$, such that 
$0$ is an  eigenvalue of the Friedrichs extension of
$L_{\lambda_1,\dots,\lambda_k,a_1,\dots,a_k}$, namely there exists $u\in
H^1(\R^N)\setminus\{0\}$ solving $L_{\lambda_1,\dots,\lambda_k,a_1,\dots,a_k} u=0$.  
\end{Theorem}

\noindent Actually, it turns out that, under assumptions  (\ref{eq:19}) and
(\ref{eq:72}), the set of configurations
$(a_1,\dots,a_k)$ for which $0$ is an  eigenvalue of
$L_{\lambda_1,\dots,\lambda_k,a_1,\dots,a_k}$, namely for which there exists an
$L^2$-bound state with null energy, even disconnects
$\R^{Nk}\setminus \Sigma$, where
$\Sigma:=\{(a_1,\dots,a_k)\in\R^{Nk}:\ a_i=a_j\ \text{for some }i\neq j\}$.

\medskip
The paper is organized as follows. In Section \ref{sec:positivity} we
give a condition for positivity of Schr\"odinger operators with
potentials in $\mathcal V$. Section \ref{sec:bound} contains the
analysis of the asymptotic  
behavior near the poles of solutions to Schr\"odinger equations with
Hardy type potentials,  
the proofs of Lemmas \ref{l:estmu} and \ref{l:sep}, and an extension to the case of
infinitely many singularities  
on reticular structures.  In Section \ref{sec:pertinf} the
possibility of perturbing positive operators at infinity is discussed,
while 
Section~\ref{sec:scattering} is devoted to the proof of the Separation
Theorem \ref{t:scattering}. Section \ref{mainresult} contains the
proof of Theorem
\ref{t:mainresult}. 
In section \ref{sec:best-constants-hardy} we study the problem of
attainability of $\mu(V)$, proving Proposition \ref{p:attai} and
discussing the continuity of $\mu(V)$ with respect to the masses and
the location of singularities.
In Section~\ref{sec:semi-bound-essent} we prove
Theorem \ref{t:self-adjo}. Finally  Section \ref{sec:spectr-schr-oper-1} contains
a detailed description of the spectrum of Schr\"odinger operators with
potentials in $\mathcal V$ and the proofs of Proposition \ref{p:spe}
and Theorem \ref{t:bottom}.

\medskip
\noindent
{\bf Notation. } We list below some notation used throughout the
paper.\par
\begin{itemize}
\item[-]$B(a,r)$ denotes the ball $\{x\in\R^N: |x-a|<r\}$ in $\R^N$ with
center at $a$ and radius $r$.
\item[-] For any $A\subset \R^N$, $\alchi_{A}$ denotes the
characteristic function of $A$. 
\item[-] $S$ is the best constant in the Sobolev inequality
$S\|u\|_{L^{2^*}(\R^N)}^2\leq \|u\|_{{\mathcal D}^{1,2}(\R^N)}^2$.
\item[-] For all $t\in\R$,  $t^+:=\max\{t,0\}$ (respectively
$t^-:=\max\{-t,0\}$) denotes the positive (respectively negative)
part of $t$.
\item[-] For all functions $f:\ \R^N\to\R$, $\mathop{\rm supp}f$
denotes the support of $f$, i.e. the closure of the set of points
where $f$ is non zero.
\item[-] $\omega_N$ denotes the volume of the unit ball in $\R^N$.
\item[-] For any open set $\Omega\subset\R^N$, ${\mathcal D}'(\Omega)$
denotes the space of distributions in $\Omega$.
\end{itemize}

\medskip
\noindent {\bf Acknowledgements. } The authors would like to thank
Maria J. Esteban and Diego Noja for fruitful discussions and helpful
suggestions. They are also grateful to Yehuda Pinchover for drawing
their attention to important references related to the present paper.

\section{A positivity condition in the class ${\mathcal V}$}\label{sec:positivity}

Thanks to Sobolev's inequality, for a Schr\"odinger operator $-\Delta -V$, $V\in
L^{N/2}(\R^N)$, a general positivity  condition is
$$
\|V^+\|_{L^{N/2}(\R^N)}<S.
$$
We mention that criteria for a potential energy operator $V$ ($V$ possibly
changing sign or even being a complex-valued distribution) to be
relative form-bounded with respect to the laplacian, i.e. satisfying
\begin{equation}\label{eq:22}
\bigg|\int_{\R^N}V(x)|u(x)|^2\,dx\bigg|\leq{\rm const}\,\int_{\R^N}|\n
u(x)|^2\,dx,
\end{equation}
are discussed in \cite{mazia}. In particular (\ref{eq:22}) implies
boundedness from below of the associated quadratic form in $\Di$  and 
positivity of Schr\"odinger operators with small multiple of $V$ as
a potential. However, this type of result cannot answer the question
arisen in the present paper about the positivity of forms
$Q_{\lambda_1,\dots,\lambda_k,a_1,\dots,a_k}$ for given masses $\lambda_1,\dots,\lambda_k$.

Furthermore, as $\Di$ is embedded in the Lorenz space $L^{2^*,2}(\R^N)$, the 
positivity of the quadratic form associated to operators $-\Delta-V$, 
$V\in L^{N/2,\infty}$, would be ensured by a smallness condition 
on $\|V^+\|_{L^{N/2,\infty}(\R^N)}$, where we recall that
$$
\|f\|_{L^{N/2,\infty}(\R^N)}:=\sup_{{\substack{X\subset\R^N\\\text{\rm
measurable}}}}\frac{{\displaystyle{\int_X|f(x)|\,dx
}}}{{\displaystyle{|X|^{1-\frac2N}}}}.
$$
We remark that potentials in the class ${\mathcal V}$ belong to
$L^{N/2,\infty}(\R^N)$, but their norm in $L^{N/2,\infty}(\R^N)$ is
not small. Indeed for each pole $a_i$, a direct calculation yields
$$
\bigg\|\frac{\alchi_{B(a_i,r_i)}}{|x-a_i|^2}\bigg\|_{L^{N/2,\infty}(\R^N)}\geq C(N)
$$
for some positive constant $C(N)$ independent of $r_i$. Hence, just by
increasing the number of poles, we can exhibit  potentials in ${\mathcal V}$ 
with as large norms as we want.
In the sequel, see Remark \ref{rem:ret}, we will provide an example of
potentials having infinite $L^{N/2,\infty}(\R^N)$-norm, but giving rise 
to positive quadratic forms.
Therefore to provide positivity conditions in the class is a
nontrivial issue. 

In this section we provide a criterion for establishing positivity of
Schr\"odinger operators with potentials in $\mathcal V$ in the spirit
of the well-known Allegretto-Piepenbrink theory
\cite{allegretto,piepenbrink}, which relates the existence of positive
solutions to a Schr\"odinger equation with the positivity of the
spectrum of the corresponding operator. For analogous criteria for
potentials in the Kato class we refer to \cite[Theorem 2.12]{CFKS}.

\begin{Lemma}\label{l:positivity_condition}
  Let $V=\sum_{i=1}^k\frac{\lambda_i\alchi_{B(a_i,r_i)}(x)}{|x-a_i|^2}+
  \frac{\lambda_{\infty} \alchi_{\R^N\setminus
      B(0,R)}(x)}{|x|^2}+W(x)\in{\mathcal V}$. Then the two following
  conditions are equivalent:
\begin{align*}
(i)\quad &\mu(V):=\inf_{u\in\Di\setminus\{0\}}\frac{\int_{\R^N}\big(|\n u(x)|^2-V(x)u^2(x)\big)\,dx}{\int_{\R^N}|\n u(x)|^2\,dx}>0;\\[10pt]
(ii)\quad &\text{there exist }\e>0 \text{ and }\varphi\in \Di,\
\varphi>0\text{ in }\R^N\setminus\{a_1,\dots,a_k\},\text{ and }\varphi\text{ 
$C^1$-smooth }\\
&\text{in }\R^N\setminus\{a_1,\dots,a_k\},\text{ such that
}-\Delta\varphi(x)-V(x)\varphi(x)>\e\, 
V^+(x)\varphi(x)\text{ a.e. in }\R^N.
\end{align*}
\end{Lemma}

\begin{pf}
We first observe that, from Hardy's, H\"older's, and Sobolev's inequalities,  
\begin{equation}\label{eq:4}
\int_{\R^N}V(x)u^2(x)\,dx\leq\bigg[\frac4{(N-2)^2} \Big(\sum_{i=1}^k\lambda_i^++\lambda_{\infty}^+\Big)
+S^{-1}\|W^+\|_{L^{N/2}(\R^N)}\bigg]\int_{\R^N}|\n u(x)|^2\,dx,
\end{equation}
for every $u\in\Di$.

Assume that $(i)$ holds. If $0<\e<\frac{\mu(V)}{2}\Big[\frac4{(N-2)^2} \Big(\sum_{i=1}^k\lambda_i^++\lambda_{\infty}^+\Big) +S^{-1}\|W^+\|_{L^{N/2}(\R^N)}\Big]^{-1}$, from (\ref{eq:4}) it follows that 
$$
\int_{\R^N}\big(|\n u(x)|^2-V(x)u^2(x)-\e\,V^+(x)u^2(x)\big)\,dx\geq \frac{\mu(V)}2\int_{\R^N}|\n u(x)|^2\,dx.
$$
As a consequence, for any fixed $p\in L^{N/2}(\R^N)\cap L^{\infty}(\R^N)$, $p(x)>0$ a.e. in $\R^N$, the infimum
$$
\nu_p(V+\e V^+)=\inf_{u\in\Di\setminus\{0\}}\dfrac{{{\int_{\R^N}\big(|\n  u(x)|^2-V(x)\,
u^2(x)-\e V^+(x)\,u(x)^2\big)\,dx}}}{{{\int_{\R^N}p(x)u^2(x)}}}
$$
is strictly positive and  attained by some function $\varphi\in\Di\setminus\{0\}$ satisfying
$$
-\D\varphi(x)-V(x)\varphi(x)=\e\,V^+(x)\varphi(x)+\nu_p(V+\e V^+)p(x)\varphi(x).
$$
By evenness we can assume $\varphi\geq0$. Since $V\in{\mathcal V}$, the
Strong Maximum Principle allows us to conclude that $\varphi>0$ in
$\R^N\setminus\{a_1,\dots,a_k\}$, while standard regularity theory
ensures regularity of $\varphi$ outside the poles. Hence $(ii)$
holds. 

Assume now that $(ii)$ holds. For any $u\in C^{\infty}_{\rm
c}(\R^N\setminus\{a_1,\dots,a_k\})$, testing the inequality
satisfied by $\varphi$ with $u^2/\varphi$ we get  
\begin{align*}
\e\int_{\R^N}V^+(x)u^2(x)\,dx
&\leq2\int_{\R^N}\frac{u(x)}{\varphi(x)}\n u(x)\cdot\n
\varphi(x)\,dx-\int_{\R^N}\frac{u^2(x)}{\varphi^2(x)}|\n
\varphi(x)|^2\,dx-\int_{\R^N}V(x)u^2(x)\,dx\\ 
&\leq\int_{\R^N}|\n u(x)|^2\,dx-\int_{\R^N}V(x)u^2(x)\,dx.
\end{align*}
By density of $C^{\infty}_{\rm c}(\R^N\setminus\{a_1,\dots,a_k\})$ in
$\Di$, we deduce that 
\begin{equation}\label{eq:3}
\inf_{u\in\Di\setminus\{0\}}\frac{\int_{\R^N}\big(|\n
u(x)|^2-V(x)u^2(x)\big)\,dx}{\int_{\R^N}V^+(x)u^2(x)\,dx}\geq\e.
\end{equation}
From \eqref{eq:3} we obtain that $\mu(V)\geq\e/(1+\e)$. Indeed, 
for every $u\in\Di\setminus\{0\}$, 
$$\int_{\R^N}\big(|\n
u(x)|^2-V(x)u^2(x)\big)\,dx\geq\e\int_{\R^N}V^+(x)u^2(x)\,dx
\geq\e\int_{\R^N}V(x)u^2(x)\,dx,$$ 
so that 
$$\int_{\R^N}V(x)u^2(x)\leq\frac{1}{1+\e}\int_{\R^N}|\n u(x)|^2\,dx,$$
implying
\begin{equation}\label{eq:66}
\int_{\R^N}\big(|\n
u(x)|^2-V(x)u^2(x)\big)\,dx\geq\frac{\e}{1+\e}\int_{\R^N}|\n u(x)|^2\,dx.
\end{equation}
\end{pf}

\section{The Shattering Lemma}\label{sec:bound}

A starting point for the study of positivity of Schr\"odinger operators with
potentials in ${\mathcal V}$, which will be also a key ingredient for
the study of their spectral structure,  relies in the Shattering Lemma~\ref{l:sep}, 
which ensures the positivity of Schr\"odinger operators
with potentials whose singularities are localized in a small
neighborhood of the corresponding poles, thus avoiding mutual interactions. 

The proof of the Shattering Lemma consists in constructing
supersolutions for each operator ${\mathcal L}^{\delta}_{\lambda_i,a_i}$
(see Lemma \ref{l:positivity_condition}) and then summing up to obtain
a supersolution for ${\mathcal
  L}^{\delta}_{\lambda_1,\dots,\lambda_k,a_1,\dots,a_k}$.  In order to take
account of the interactions, we need to evaluate the exact behavior of
such supersolutions at each pole and at $\infty$, as described
in Lemmas \ref{l:1} and \ref{l:2} below, for the proof of which we
refer to \cite{FS3} and \cite{murata}, see also \cite{pinchover94}.

\begin{Lemma}\label{l:1}
Let $\varphi\in\Di$, $\varphi\geq 0$ a.e. in $\R^N$,
$\varphi\not\equiv 0$, be a weak solution of  
\begin{equation}\label{eq:50}
-\Delta\varphi(x)=\bigg[\frac{\lambda\,\alchi_{B(0,1)}(x)}{|x|^2}+
q(x)\bigg]\,\varphi(x)\quad\text{ in }\R^N,
\end{equation}
where $\lambda<\frac{(N-2)^2}{4}$ and
$q\in L^{\infty}_{\rm loc}(\R^N\setminus\{0\})$. Then 
\begin{itemize}
\item[$(i)$] if $q(x)=O(|x|^{-(2-\e)})$ as $|x|\to0$ for some $\e>0$,
  there exists a positive constant $C$ (depending on $q$, $\lambda$,
  $\e$, and $\varphi$) such that
\begin{align*}
\frac1C|x|^{-a_{\lambda}}\leq \varphi(x)\leq
C|x|^{-a_{\lambda}}\quad\text{ for all }x\in B(0,1),
\end{align*}
where $a_{\lambda}=\frac{N-2}2-\sqrt{\big(\frac{N-2}2\big)^2-\lambda}$;
\item[$(ii)$] if $q(x)=O(|x|^{-2-\e})$ as $|x|\to +\infty$ for some $\e>0$, there exists a positive constant $C$ (depending on $q$, $\lambda$, $\e$, and $\varphi$) such that 
\begin{align*}
\frac1C|x|^{-(N-2)}\leq \varphi(x)\leq
C|x|^{-(N-2)}\quad \text{ for all }x\in\R^N\setminus B(0,1).
\end{align*}
\end{itemize}
\end{Lemma}

\begin{remark}\label{r:h1}
Scanning through the proof of the regularity result contained in \cite[Theorem
1.1]{FS3}, it is possible to
clarify the dependence of the estimates stated  in the above
lemma on the data of the
problem. Indeed it turns out that if $\varphi$ solves (\ref{eq:50})
and $u$ is given by  
\begin{equation*}
u(x)=|x|^{a_{\lambda}}\varphi(x),
\end{equation*}
 then, by \cite[Theorem
1.1]{FS3}, $u\in C^{0,\a}(B(0,1))$ for some $\alpha>0$ and  $\|u\|_{C^{0,\a}(B(0,1))}\leq
c\,\|\varphi\|_{H^1(B(0,R))}$ for any $R>1$ and for some positive constant $c$ (depending
on $R$, $q$, $\e$, $N$ and $\lambda$) which stays bounded uniformly with respect to $\lambda$
whenever $\lambda$ stays bounded from below and above away from $(N-2)^2/4$. Hence 
$$
\varphi(x)\leq  c(\lambda,N,q,\e,R)\,|x|^{-a_{\lambda}}\|\varphi\|_{H^1(B(0,R))}.
$$
In particular, we can bound $c(\lambda,N,q,\e,R)$ uniformly with respect to
$\lambda$ if $\lambda$ varies in a compact subset of $(-\infty,(N-2)^2/4)$.
\end{remark}

\begin{Lemma}\label{l:2}
  Let $\varphi\in\Di$, $\varphi\geq 0$ a.e. in $\R^N$,
  $\varphi\not\equiv0$, be a weak solution of
\begin{equation*}
-\Delta\varphi(x)=\bigg[ \frac{\lambda\,
\alchi_{\R^N\setminus B(0,1)}(x)}{|x|^2}+
h(x)
\bigg]\, \varphi(x)\quad\text{ in
}\R^N, 
\end{equation*}
where  $\lambda<\frac{(N-2)^2}{4}$ and 
$h\in L^{\infty}_{\rm loc}(\R^N\setminus\{0\})$. Then 
\begin{itemize}
\item[$(i)$] if $h(x)=O(|x|^{-2-\e)})$ as $|x|\to+\infty $ for some
  $\e>0$, there exists a positive
constant $C$ (depending on $h$, $\lambda$, $\e$, and $\varphi$) such that 
\begin{align*}
\frac1C\, |x|^{-(N-2-a_{\lambda})}\leq\varphi(x)\leq C\,
|x|^{-(N-2-a_{\lambda})}\quad\text{for all }x\in\R^N\setminus
B(0,1),
\end{align*}
where $a_{\lambda}=\frac{N-2}2-\sqrt{\big(\frac{N-2}2\big)^2-\lambda}$;
\item[$(ii)$] if $h(x)=O(|x|^{-(2-\e)})$ as $|x|\to0$ for some
  $\e>0$, there exists a positive
constant $C$ (depending on $h$, $\lambda$, $\e$, and $\varphi$) such that 
\begin{align*}
\frac{1}{C}\leq\varphi(x)\leq C
\quad\text{for all }x\in B(0,1).
\end{align*}
\end{itemize}
\end{Lemma}

\noindent 
We now prove the bound from above of the first
$\Di$-eigenvalue of Schr\"odinger operators with potential 
in ${\mathcal V}$ stated in Lemma \ref{l:estmu}.

\bigskip\noindent

\begin{pfn}{Lemma \ref{l:estmu}}
Let us first consider the case in which $\lambda_i<0$ for all
$i=1,\dots,k,\infty$. Let us fix
$u\in C^{\infty}_{\rm c}(\R^N\setminus\{a_1,\dots,a_k\})$ and
$P\in B(0,R_0)\setminus\{a_1,\dots,a_k\}$. Letting
$u_{\mu}(x)=\mu^{-\frac{N-2}2}u\big(\frac{x-P}\mu\big)$, for 
$\mu$ small there holds 
$$
\mu(V)\leq 1-\frac{\int_{\R^N}W(x)u_{\mu}^2(x)\,dx}
{\int_{\R^N}|\n u_{\mu}(x)|^2\,dx}=1+o(1)\quad\text{as }\mu\to 0^+.
$$
Letting $\mu\to 0^+$ we obtain that $\mu(V)\leq 1$.

Assume now $\max_{i=1,\dots,k,\infty}\lambda_i>0$. 
Suppose $\lambda_1\leq\lambda_2\leq\dots\lambda_k$ and let $\e>0$. From optimality of
the best constant in the classical Hardy inequality (\ref{eq:hardy})
and by density of $C^{\infty}_{\rm c}(\R^N\setminus\{a_1,\dots,a_k\})$
in $\Di$, there exists $\phi\in C^{\infty}_{\rm
  c}(\R^N\setminus\{a_1,\dots,a_k\})$ such that  
\begin{equation*}
\int_{\R^N}|\n\phi(x)|^2dx<\bigg[\frac{(N-2)^2}4+\e\bigg]\int_{\R^N}\frac{\phi^2(x)}{|x|^2}\,dx.
\end{equation*}
Letting $\phi_{\mu}(x)=\mu^{-\frac{N-2}2}\phi\big(\frac{x-a_k}\mu\big)$, for any $\mu>0$ there holds
\begin{align*}
\mu(V)&\leq 1-\lambda_k\frac{\int_{B(a_k,r_k)}|x-a_k|^{-2}\phi_{\mu}^2(x)\,dx}{\int_{\R^N}|\n \phi_{\mu}(x)|^2\,dx}
-\sum_{i\neq k}\lambda_i
\frac{\int_{B(a_i,r_i)}|x-a_i|^{-2}\phi_{\mu}^2(x)\,dx}{\int_{\R^N}|\n \phi_{\mu}(x)|^2\,dx}\\
&\quad-
\lambda_{\infty}\frac{\int_{\R^N\setminus B(0,R_0)}|x|^{-2}\phi_{\mu}^2(x)\,dx}{\int_{\R^N}|\n \phi_{\mu}(x)|^2\,dx}-\frac{\int_{\R^N}W(x)\phi_{\mu}^2(x)\,dx}
{\int_{\R^N}|\n \phi_{\mu}(x)|^2\,dx}
\\
&=1-\lambda_k\frac{\int_{\R^N}|x|^{-2}\phi^2(x)\,dx}{\int_{\R^N}|\n \phi(x)|^2\,dx}+o(1)\quad\text{as}\quad \mu\to 0^+.
\end{align*}
Letting $\mu\to 0^+$, by the choice of $\phi$ we obtain
$$
\mu(V)\leq 1-\lambda_k\bigg[\frac{(N-2)^2}4+\e\bigg]^{-1}
$$
for any $\e>0$. Letting $\e\to 0$ we derive that
$\mu(V)\leq 1-\frac{4\lambda_k}{(N-2)^2}$. Repeating the same argument with 
$\phi_{\mu}(x)=\mu^{-\frac{N-2}2}\phi(x/\mu)$ and letting $\mu\to+\infty$ we obtain also that $\mu(V)\leq 1-\frac{4\lambda_{\infty}}{(N-2)^2}$. The required estimate is thereby proved.
\end{pfn}

\begin{pfn}{Lemma \ref{l:sep}}
Assume that $\max_{i=1,\dots,k,\infty}\lambda_i>0$, otherwise there is nothing to
prove. Let us fix $0<\e<\frac{(N-2)^2}{4\max_{i=1,\dots,k,\infty}\lambda_i}-1$, so
that 
$$
\widetilde\lambda_{\infty}:=\lambda_{\infty}+\e\lambda_{\infty}^+<\frac{(N-2)^2}4\quad\text{and}\quad\widetilde\lambda_i:=\lambda_{i}+\e\lambda_{i}^+<\frac{(N-2)^2}4\quad\text{ for all }i=1,\dots,k.
$$
By scaling properties of the operator and in view of Lemma
\ref{l:positivity_condition}, to prove the statement it is enough to find $\varphi\in\Di$ positive and smooth outside the singularities such that 
\begin{equation}\label{eq:80}
-\Delta\varphi(x)-\sum_{i=1}^kV_i(x)\,\varphi(x)- V_{\infty}(x)\,\varphi(x)>0\quad\text{a.e. in }\R^N,
\end{equation}
where
$$
V_i(x)=\frac{\widetilde\lambda_i\,\alchi_{B(a_i/\delta,1)}(x)}
{\big|x-\frac{a_i}{\delta}\big|^2},\quad
V_{\infty}(x)=\frac{\widetilde\lambda_{\infty}\,\alchi_{\R^N\setminus
    B(0,R_0/\delta)}(x)}{|x|^2}, 
$$ 
and $\d>0$ depends only on the location of poles and  on $\e$. 
Indeed, if (\ref{eq:80}) holds for some positive $\varphi$, then Lemmas
\ref{l:positivity_condition} and \ref{l:estmu} ensure that 
 $$
1-\frac4{(N-2)^2}\,\max_{i=1,\dots,k,\infty}\lambda_i
\geq \mu\bigg(\sum_{i=1}^k\frac{\lambda_i\,\alchi_{B(a_i,\delta)}(x)}{|x-a_i|^{2}}+ 
 \frac{\lambda_{\infty}\alchi_{\R^N\setminus
     B(0,R_0)}(x)}{|x|^{2}}\bigg)\geq\frac{\e}{1+\e}
 $$
for all $0<\e<\frac{(N-2)^2}{4\max_{i=1,\dots,k,\infty}\lambda_i}-1$ and
the result follows letting $\e\to
\frac{(N-2)^2}{4\max_{i=1,\dots,k,\infty}\lambda_i}-1$.  

In order to find a positive supersolution to (\ref{eq:50}), for all
$i=1,\dots,k$ let us set, for some $0<\tau<1$,
$$
p_i(x):=p\Big(x-\frac{a_i}{\delta}\Big)\quad\text{ where }p(x)=\frac1{|x|^{2-\tau}(1+|x|^2)^{\tau}},
$$
and
$$
p_{\infty}(x)=\frac{\d^{\tau}R_0^{-\tau}}{|x|^{2-\tau}(1+|\delta x/R_0|^2)^{\tau}}.
$$
Since $p_i,p_{\infty}\in L^{N/2}(\R^N)$, it is easy to see that the weighted first eigenvalue
$$
\mu_i=\min_{\Di\setminus\{0\}}\dfrac{{\displaystyle{\int_{\R^N}\big(|\n
\varphi(x)|^2-V_i(x)\,\varphi^2(x)}}\big)\,dx}{{\displaystyle{\int_{\R^N}p_i(x)\varphi^2(x)}}\,dx}
$$
is positive and attained by some function $\varphi_i\in\Di$,
$\varphi_i>0$ and smooth  in $\R^N\setminus\{a_i/\d\}$ and also
\begin{align*}
\mu_{\infty}&=\min_{\Di\setminus\{0\}}\frac{{\displaystyle{\int_{\R^N}\big(|\n
\varphi(x)|^2-V_{\infty}(x)\,
\varphi^2(x)}}\big)\,dx}{{\displaystyle{\int_{\R^N}p_{\infty}(x)\varphi^2(x)}}\,dx}\\
&=\min_{\Di\setminus\{0\}}\frac{{\displaystyle{\int_{\R^N}|\n
\varphi(x)|^2\,dx-\tilde \lambda_{\infty}\int_{\R^N\setminus B(0,1)}
\frac{\varphi^2(x)}{|x|^2}}\,dx}}{{\displaystyle{\int_{\R^N}p(x)\varphi^2(x)}}\,dx}
\end{align*}
is positive and attained by some function $\varphi_{\infty}\in\Di$,
$\varphi_{\infty}>0$  and smooth  in $\R^N\setminus\{0\}$.
The function $\varphi_i$ satisfy
$$
-\D\varphi_i(x)-V_i(x)\,\varphi_i(x)=\mu_i p_i(x)\varphi_i(x)
$$
while $\varphi_{\infty}$ satisfy
$$
-\D\varphi_{\infty}(x)-V_{\infty}(x)\,\varphi_{\infty}(x)=\mu_{\infty}
p_{\infty}(x)\varphi_{\infty}(x),
$$
hence $\psi(x):=\varphi_{\infty}\big(\frac{R_0}{\delta}x\big)$ satisfies
$$
-\Delta\psi(x)- \frac{\tilde\lambda_{\infty}\,
\alchi_{\R^N\setminus B(0,1)}(x)}{|x|^2}
\, \psi(x)=\mu_{\infty}
p(x)\psi(x).
$$
Lemmas \ref{l:1} and \ref{l:2} yield a constant $C_0>0$ (independent on
$\d$) such that 
\begin{align}\label{eq:32}
\frac1{C_0}\Big|x-\frac{a_i}\d\Big|^{-a_{\tilde\lambda_i}}\leq \varphi_i(x)\leq
C_0\Big|x-\frac{a_i}{\d}\Big|^{-a_{\tilde\lambda_i}},\hskip2.05cm&\text{ for
  all }x\in 
B(a_i/\d,1)\\[5pt]
\label{eq:34}
\frac1{C_0}\Big|x-\frac{a_i}{\d}\Big|^{-(N-2)}\leq \varphi_i(x)\leq
C_0\Big|x-\frac{a_i}{\d}\Big|^{-(N-2)},\hskip.9cm\quad&\text{ for all
}x\in\R^N\setminus B(a_i/\d,1) ,\\[5pt]
\label{eq:35}\frac1{C_0}\,
\Big|\frac{\d
  x}{R_0}\Big|^{-(N-2-a_{\tilde\lambda_{\infty}})}\leq\varphi_{\infty}(x)\leq
C_0\, \Big|\frac{\d x}{R_0}\Big|^{-(N-2-a_{\tilde\lambda_{\infty}})},\quad&\text{ for all
}x\in\R^N\setminus 
B(0,R_0/\d),\\[5pt]  
\label{eq:36}
\frac{1}{C_0}\leq\varphi_{\infty}(x)\leq C_0,\hskip4.3cm
\hskip1cm\quad&\text{ for all }x\in B(0,R_0/\d).
\end{align}
Let $\varphi=\sum_{i=1}^k\varphi_i+\eta\varphi_{\infty}$ for some
$0<\eta<\inf\Big\{
\frac{\strut \mu_i}{\strut 4C_0^2\tilde\lambda_i}
:\  i=1,\dots,k,\ \tilde\lambda_i>0\Big\}$. Then we have 
\begin{align*}
-\D\varphi(x)-\sum_{i=1}^kV_i(x)\,\varphi(x)-
 V_{\infty}(x)\,\varphi(x)&=\sum_{i=1}^k\mu_ip_i(x)\varphi_i(x)+
\mu_{\infty}p_{\infty}(x)\eta\,\varphi_{\infty}(x)\\
&-\sum_{i\neq j}V_i(x)\varphi_j(x)-\eta\sum_{i=1}^kV_i(x)\varphi_{\infty}(x)-V_{\infty}(x)
\sum_{i=1}^k\varphi_i(x).
\end{align*}
In particular a.e. in the set $B(0,R_0/\d)\setminus\bigcup_{i=1}^k B(a_i/\d,1)$, we have
$$
-\D\varphi(x)-\sum_{i=1}^kV_i(x)\,\varphi(x)- V_{\infty}(x)\,\varphi(x)=\sum_{i=1}^k\mu_ip_i(x)\varphi_i(x)+
\mu_{\infty}p_{\infty}(x)\eta\,\varphi_{\infty}(x)>0.
$$
Let us consider $B(a_i/\d,1)$. 
If $\widetilde\lambda_i<0$, we have easily that 
$$
-\D\varphi-\sum_{i=1}^kV_i(x)\,\varphi(x)-V_{\infty}(x)\,\varphi(x)>0\quad\mbox{ a.e. in $B(a_i/\d,1)$}.
$$
If $\widetilde\lambda_i>0$ we can choose $\tau<a_{\widetilde\lambda_i}$.
Since, for $\delta$ small, 
$$B(a_i/\d,1)\subset B(0,R_0/\d)\quad\text{and}\quad
B(a_i/\d,1)\subset\R^N\setminus B(a_j/\d,1)\text{ for }j\neq i,
$$
from
(\ref{eq:32}), (\ref{eq:34}), and (\ref{eq:36})  it
follows that,
in  $B(a_i/\d,1)$, 
\begin{align*}
-\D\varphi(x)-&\sum_{i=1}^kV_i(x)\,\varphi(x)-V_{\infty}(x)\,\varphi(x)\geq
\mu_ip_i(x)\varphi_i(x)-V_i(x)\Big(\sum_{j\neq
  i}\varphi_j(x)+\eta\,\varphi_{\infty}(x)\Big)\\ 
&\geq \Big|x-\frac{a_i}{\delta}\Big|^{-2}
\bigg[\frac{\mu_i}{2^{\tau}C_0}\Big|x-\frac{a_i}{\delta}\Big|^{\tau-a_{\widetilde\lambda_i}}
-C_0\tilde\lambda_i\bigg(\sum_{j\neq i}\Big|x-\frac{a_j}{\delta}\Big|^{-(N-2)}+\eta\bigg)\bigg].
\end{align*}
It is easy to see that for $\delta$ small 
$$
\Big|x-\frac{a_i}{\delta}\Big|^{\tau-a_{\widetilde\lambda_i}}\geq
1\quad\text{and}\quad
\Big|x-\frac{a_j}{\delta}\Big|^{-(N-2)}\leq\Big(\frac{2}{|a_i-a_j|}\Big)^{N-2}
\d^{N-2}<\frac{\eta}{k-1}, 
$$
and hence the choice of $\eta$ ensures that
$$
-\D\varphi(x)-\sum_{i=1}^kV_i(x)\,\varphi(x)-V_{\infty}(x)\, \varphi(x)>0\quad\text{ a.e. in }B(a_i/\d,1),
$$
provided $\delta$ is sufficiently small.

Let us finally consider $\R^N\setminus B(0,R_0/\d)$. 
If $\widetilde\lambda_{\infty}<0$, then
$$
-\D\varphi(x)-\sum_{i=1}^kV_i(x)\,\varphi(x)-V_{\infty}(x)\,\varphi(x)>0
\quad\text{a.e. in }\R^N\setminus B(0,R_0/\d).
$$
If $\widetilde\lambda_{\infty}>0$ we can choose $\tau<a_{\tilde\lambda_{\infty}}$.
From (\ref{eq:34}--\ref{eq:35}), we deduce that in $\R^N\setminus
B(0,R_0/\d)$ 
\begin{align*}
-\D\varphi(x)-&\sum_{i=1}^kV_i(x)\,\varphi(x)-V_{\infty}(x)\,\varphi(x)\geq
\mu_{\infty}p_{\infty}(x)\eta\,\varphi_{\infty}(x)- V_{\infty}(x)\sum_{i=1}^k\varphi_i(x)\\
&\geq \frac1{|x|^{2}} \bigg[
\frac{\mu_{\infty}\eta}{2^{\tau}C_0}\Big|\frac{\d
  x}{R_0}\Big|^{-(N-2-a_{\tilde\lambda_{\infty}})-\tau}-C_0\tilde \lambda_{\infty}\sum_{i=1}^k \Big|x-\frac{a_i}{\delta}\Big|^{-(N-2)}\bigg]. 
\end{align*}
It is easy to see that,  in
$\R^N\setminus B(0,R_0/\d)$, 
$\big|x-\frac{a_i}{\delta}\big|\geq\big(1-\frac{\a}{R_0}\big)|x|$
where $\a=\max\{|a_j|\}_j$, hence 
\begin{align*}
 -\D\varphi(x)-&\sum_{i=1}^kV_i(x)\,\varphi(x)-V_{\infty}(x)\,\varphi(x)\geq
 \mu_{\infty}p_{\infty}(x)\eta\,\varphi_{\infty}(x)- V_{\infty}(x)\sum_{i=1}^k\varphi_i(x)\\
 &\geq \frac1{|x|^{N}} \bigg[
 \frac{\mu_{\infty}\eta}{2^{\tau}C_0}\Big|\frac{\d}{R_0}
\Big|^{-(N-2)} 
-C_0\tilde \lambda_{\infty}k\Big(1-\frac{\a}{R_0}\Big)^{-(N-2)}\bigg]>0  \quad\text{ a.e. in }\R^N\setminus
B(0,R_0/\d)
\end{align*}
provided $\delta$ is sufficiently small.
The proof is thereby complete.
\end{pfn}

\begin{remark}\label{r:shattering-lemma}
  For $(\lambda_1,\dots,\lambda_k,\lambda_{\infty})\in
  \big(-\infty,{(N-2)^2}/4\big)^{k+1}$, let
  $\{a_1^n,a_2^n,\dots,a_k^n\}\subset B(0,R_0)\subset\R^N$ be a
  sequence of configurations approximating
  $\{a_1,a_2,\dots,a_k\}\subset B(0,R_0)\subset\R^N$, $a_i\neq a_j$
  for $i\neq j$, i.e. $a_i^n\to a^i$ as $n\to\infty$ for all
  $i=1,\dots,k$.  From the proofs of Lemmas
  \ref{l:positivity_condition} and \ref{l:sep}, for any
  $0<\alpha<1-\frac{4\max_{i=1,\dots,k,\infty}\lambda_i}{(N-2)^2}$, we can
  choose $\delta>0$ independently of $n$ such that
\begin{equation*}
\mu(\widetilde V_n)\geq
\begin{cases}
1-\frac{4\max_{i=1,\dots,k,\infty}\lambda_i}{(N-2)^2}-\alpha,&\text{if
}\max_{i=1,\dots,k,\infty}\lambda_i>0,\\
1,&\text{if }\max_{i=1,\dots,k,\infty}\lambda_i\leq0,
\end{cases}
\end{equation*}
where
$$
\widetilde
V_n(x)=\sum_{i=1}^k\frac{\lambda_i\alchi_{B(a_i^n,\d)}}{|x-a_i^n|^2}+
\frac{\lambda_{\infty}\alchi_{\R^N\setminus  B(0,R_0)}}{|x|^2}. 
$$
\end{remark}

\noindent Let us now deal with the case of infinitely many
singularities distributed on
reticular structures.

\begin{Lemma}{\bf [\,Shattering of reticular singularities\,]}\label{l:ret}
Let $\lambda<(N-2)^2/4$ and let $\{a_n\}_n\subset\R^N$ satisfy 
\begin{equation}\label{eq:55}
\sum_{n=1}^{\infty}|a_n|^{-(N-2)}<+\infty,\quad 
\sum_{k=1}^{\infty}|a_{n+k}-a_n|^{-(N-2)}\quad\text{is bounded uniformly
  in }n,
\end{equation}
and $|a_n-a_m|\geq 1$
for all $n\neq m$. Then  there exists $\delta>0$ 
such that 
$$
\inf_{{\substack{u\in \Di
\\u\not\equiv0}}}\frac{\displaystyle\int_{\R^N}\big(|\n
u(x)|^2-V(x)u^2(x)\big)\,dx}{\displaystyle\int_{\R^N}|\n
u(x)|^2\,dx}>0
$$
where
$$
V(x)=\lambda\sum_{n=1}^{\infty}\frac{\alchi_{B(a_n,\delta)}(x)}{|x-a_n|^2}.
$$
\end{Lemma}
\begin{pf}
Let $\e>0$ such that $\tilde\lambda=\lambda+\e<(N-2)^2/4$. Arguing as in the proof of Lemma
\ref{l:sep}, we can construct  
 a function $\psi\in\Di$,
$\psi>0$ and smooth  in $\R^N\setminus\{0\}$ such that 
$$
-\D\psi(x)-\tilde\lambda\,\frac{\alchi_{B(0,1)}(x)}
{|x|^2}\,\psi(x)=\mu\,p(x)\psi(x), 
$$ 
where $\mu>0$ and $p(x)={|x|^{\tau-2}(1+|x|^2)^{-\tau}}$, for some $0<\tau<1$. Moreover, by Lemmma \ref{l:1},
\begin{align*}
\frac{|x|^{-a_{\tilde\lambda}}}C\leq \psi(x)\leq C|x|^{-a_{\tilde\lambda}}\quad\text{in }B(0,1),\quad
\text{and}\quad\frac{|x|^{-(N-2)}}C\leq \psi(x)\leq
C|x|^{-(N-2)}\quad\text{in }\R^N\setminus B(0,1).
\end{align*}
Let
$\varphi(x)=\sum_{n=1}^{\infty}\psi\big(x-\frac{a_n}\delta\big)$. For
any compact set $K$, we have that there exists $\bar n$ such that, for
all $n\geq \bar n$ and $x\in K$, $\big|\psi\big(x-\frac{a_n}\delta\big)\big|\leq
C\big|x-\frac{a_n}\delta\big|^{-(N-2)}\leq{\rm const\,}\big|\frac{a_n}\delta\big|^{-(N-2)}$. 
Then 
$$
\varphi\big|_K(x)=\sum_{n=1}^{\bar
  n-1}\psi\Big(x-\frac{a_n}\delta\Big)+\sum_{n=\bar
  n}^{\infty}\psi\Big(x-\frac{a_n}\delta\Big)\in H^1(K)+L^{\infty}(K).
$$
In particular $\varphi\in L^1_{\rm loc}(\R^N)$. Moreover 
\begin{multline*}
-\Delta\varphi(x)-\tilde\lambda\,\sum_{n=1}^{\infty}
\frac{\alchi_{B(a_n/\delta,1)}(x)}{|x-\frac{a_n}\delta|^2}
\varphi(x)\\
=\mu\sum_{n=1}^{\infty}p\Big(x-\frac{a_n}\delta\Big)
\psi\Big(x-\frac{a_n}\delta\Big)
-\tilde \lambda \sum_{m\neq
  n}\frac{\alchi_{B(a_m/\delta,1)}(x)}{|x-\frac{a_m}\delta|^2}
\psi\Big(x-\frac{a_n}\delta\Big).
\end{multline*}
Then a.e. in the set $\R^N\setminus\bigcup_{n=1}^{\infty} B(a_n/\d,1)$, we have
$$
-\Delta\varphi(x)-\tilde\lambda\,\sum_{n=1}^{\infty}
\frac{\alchi_{B(a_n/\delta,1)}(x)}{|x-\frac{a_n}\delta|^2}
\varphi(x)>0.
$$
Assume $\tilde \lambda>0$, otherwise there is nothing to prove, thus we can choose 
$\tau<a_{\tilde\lambda}$. Therefore
in each ball $B(a_n/\delta,1)$ 
\begin{align*}
-\Delta\varphi(x)&-\tilde\lambda\,\sum_{n=1}^{\infty}
\frac{\alchi_{B(a_n/\delta,1)}(x)}{|x-\frac{a_n}\delta|^2}
\varphi(x)\\
&\geq \mu \,p\Big(x-\frac{a_n}\delta\Big)
\psi\Big(x-\frac{a_n}\delta\Big)-\tilde \lambda
\frac{\alchi_{B(a_n/\delta,1)}(x)}{|x-\frac{a_n}\delta|^2} 
\sum_{m\neq
  n}\psi\Big(x-\frac{a_m}\delta\Big)\\
&\geq {\rm const\,}\Big|x-\frac{a_n}\delta\Big|^{-2}\bigg({\rm const\,}
\Big|x-\frac{a_n}\delta\Big|^{-a_{\tilde \lambda}+\e}-\sum_{m\neq
  n}\Big|x-\frac{a_m}\delta\Big|^{-(N-2)}\bigg).
\end{align*}
Since, for small $\d$, $\big|x-\frac{a_m}\delta\big|\geq
\frac{|a_m-a_n|}\d-1\geq \frac{|a_m-a_n|}{2\d}$ provided $\d$ small
enough, we deduce that
$$
\sum_{m\neq
  n}\Big|x-\frac{a_m}\delta\Big|^{-(N-2)}\leq (2\d)^{N-2}\sum_{m\neq
  n}|a_m-a_n|^{-(N-2)}\leq {\rm const\,}\d^{N-2}.
$$
Hence, we can choose $\delta$ small enough independently of $n$ such that 
\begin{align*}
-\Delta\varphi(x)&-\tilde\lambda\,\sum_{n=1}^{\infty}
\frac{\alchi_{B(a_n/\delta,1)}(x)}{|x-\frac{a_n}\delta|^2}
\varphi(x)>0
\end{align*}
a.e. in $B(a_n/\delta,1)$.
Hence we have constructed  a supersolution  
$\varphi\in L^1_{\rm loc}(\R^N)$, $\varphi>0$ in
$\R^N\setminus\{a_n\}_{n\in\N}$ and $\varphi$ smooth in
$\R^N\setminus\{a_n\}_{n\in\N}$, such that
$$
-\Delta\varphi(x)-\tilde\lambda\,\sum_{n=1}^{\infty}
\frac{\alchi_{B(a_n/\delta,1)}(x)}{|x-\frac{a_n}\delta|^2}
\varphi(x)>0\quad\text{a.e. in }\R^N,
$$
where $\tilde\lambda=\lambda+\e<(N-2)^2/4$.
Therefore, arguing as in Lemma \ref{l:positivity_condition} and taking into
account the scaling properties of the operator, we  obtain 
$$
\tilde\mu(V):=\inf_{{\substack{u\in C^{\infty}_{\rm
        c}\left(\R^N\setminus\{a_n\}_{n\in\N}\right) 
\\u\not\equiv0}}}\frac{\displaystyle\int_{\R^N}\big(|\n
u(x)|^2-V(x)u^2(x)\big)\,dx}{\displaystyle\int_{\R^N}|\n
u(x)|^2\,dx}>0,
$$
i.e.
\begin{equation}\label{eq:56}
\int_{\R^N}V(x)u^2(x)\,dx\leq(1-\tilde\mu(V))\int_{\R^N}|\n
u(x)|^2\,dx
\end{equation}
for all $u\in C^{\infty}_{\rm
        c}\left(\R^N\setminus\{a_n\}_{n\in\N}\right)$.
By density of $C^{\infty}_{\rm
        c}\left(\R^N\setminus\{a_n\}_{n\in\N}\right)$ in $\Di$ and the
      Fatou Lemma, we can easily prove that (\ref{eq:56}) holds for
      all $u\in\Di$.

\end{pf}

\begin{remark}\label{rem:ret}
The Lemma above can be used to construct examples of potentials having infinite $L^{N/2,\infty}(\R^N)$-norm,
but giving rise to positive quadratic forms.
\end{remark}

\begin{remark}\label{rem:ret2}
If the singularities $a_n$ are located on a periodic $M$-dimensional
reticular structure, $M\leq N$, i.e. if
$$
\{a_n:\ n\in\N\}=\{(x_1,x_2,\dots,x_M,0,\dots,0):\ x_i\in\Z\ \text{for
  all }i=1,\dots,M\},
$$
then $\sum_{k=1}^{\infty}|a_{n+k}-a_n|^{-(N-2)}<+\infty$ if
and only if $\sum_{k=1}^{\infty}k^{-(N-2)+M-1}<+\infty$, i.e. for $M<N-2$.
\end{remark}

\begin{remark}\label{rem:ret3}
From Lemma \ref{l:ret}, it follows that, for $\d$ small
and any $u\in\Di$, the series
$$
\sum_{n=1}^{\infty}\int_{B(a_n,\d)}\frac{u^2(x)}{|x-a_n|^2}\,dx  
$$ 
converges and 
$$
\sum_{n=1}^{\infty}\int_{B(a_n,\d)}\frac{u^2(x)}{|x-a_n|^2}\,dx \leq 
{\rm const }\int_{\R^N}|\n u(x)|^2\,dx.
$$
\end{remark}

\section{Perturbation at infinity}\label{sec:pertinf}

In this section we discuss the stability of positivity with respect to
perturbations of the potentials with a small singularity sitting at
infinity.   

\begin{Lemma}\label{l:pertinf}
Let 
$$
V(x)=\sum_{i=1}^k\frac{\lambda_i\alchi_{B(a_i,r_i)}(x)}{|x-a_i|^2}+ \frac{\lambda_{\infty}\alchi_{\R^N\setminus B(0,R)}(x)}{|x|^2}+W(x)\in{\mathcal V},
$$
where $W\in L^{\infty}(\R^N)$, $W(x)=O(|x|^{-2-\delta})$, with
$\delta>0$, as $|x|\to\infty$. Assume that $\mu(V)>0$, namely that the quadratic form associated to the operator $-\Delta-V$ is positive definite. Let $\gamma_{\infty}\in\R$ such that
$\lambda_{\infty}+\gamma_{\infty}<\big(\frac{N-2}2\big)^2$. Then there exist $\tilde R>R$ and $\Phi\in\Di$, $\Phi>0$ in $\R^N\setminus\{a_1,\dots,a_k\}$ and smooth in $\R^N\setminus\{a_1,\dots,a_k\}$ such that 
$$
-\Delta\Phi(x)-V(x)\Phi(x)-\frac{\gamma_{\infty}}{|x|^2}\alchi_{\R^N\setminus B(0,\tilde R)}\Phi(x)>0.
$$
\end{Lemma}

\begin{pf}
Let us fix 
$$
0<\e<\min\bigg\{
\sqrt{\Big(\frac{N-2}2\Big)^2-\lambda_{\infty}},\
\sqrt{\Big(\frac{N-2}2\Big)^2-\lambda_{\infty}-\gamma_{\infty}}\bigg\},
$$
$C_0>0$  such that 
$$
W(x)\leq\frac{C_0}{|x|^{2+\delta}},\quad\text{ in }\R^N,
$$ 
and  $R_0>0$  such that 
$$
\bigcup_{i=1}^k B(a_i,r_i)\subset
B(0,R_0)\quad\text{and} 
\quad 
\left(\frac{N-2}{2}\right)^2-\e^2-\lambda_{\infty}\geq
\frac{C_0}{R_0^{\delta}}. 
$$
 Let $\varphi_1\geq0$ be a smooth function such that
$$\varphi_1\equiv0\quad\text{in } B(0,R_0),\quad
\varphi_1(x)=\frac{1}{|x|^{\frac{N-2}{2}+\e}}\quad\text{in }
\R^N\setminus B(0,2R_0).
$$
For every
 $\tilde R>2R_0$, we have that 
\begin{align}\label{eq:26}
-\Delta\varphi_1(x)-\frac{\gamma_{\infty}\alchi_{\R^N\setminus
    B(0,\tilde R)}(x)}{|x|^2}\varphi_1(x)=
&\left[\left(\frac{N-2}{2}\right)^2-\e^2\right]\frac{\alchi_{B(0,\tilde R)\setminus
    B(0,2R_0)}(x)}{|x|^2}\varphi_1(x)+f_1(x)\\ 
&\notag +\left[\left(\frac{N-2}{2}\right)^2-\e^2-
  \gamma_{\infty}\right]\frac{\alchi_{\R^N\setminus
    B(0,\tilde R)}(x)}{|x|^2}\varphi_1(x), 
\end{align}  
where $f_1$ is a smooth function with compact support. 
Let us choose a smooth function with compact support $f_2$ such that 
\begin{align*}
&f_1+f_2\geq 0\text{ in }\R^N,\quad f_1+f_2>0\text{ in }B(0,2R_0), \quad\text{and}\\
& f_2+W\alchi_{B(0,2
  R_0)}\varphi_1+\frac{\lambda_{\infty}}{|x|^2}\alchi_{B(0,2 R_0)\setminus B(0, R_0)
}\varphi_1\geq0\text{ in }\R^N.
\end{align*}
Since $f_2+W\alchi_{B(0,2 R_0)}\varphi_1+{\lambda_{\infty}}{|x|^{-2}}\alchi_{B(0,2
  R_0)\setminus B(0, R_0)}\varphi_1\in
L^{\frac{2N}{N+2}}(\R^N)$
and $\mu(V)>0$,
in view of the Lax-Milgram Theorem there exists $\varphi_2\in\Di$, $\varphi_2>0$ in
$\R^N\setminus\{a_1,\dots,a_k\}$, satisfying
\begin{equation}\label{eq:27}
-\Delta\varphi_2(x)-V(x)\varphi_2(x)=
f_2(x)+\bigg[W(x)\alchi_{B(0,2R_0)}(x)+\frac{\lambda_{\infty}\alchi_{B(0,2
  R_0)\setminus B(0, R_0)}(x)}{|x|^2}\bigg]\varphi_1(x).
\end{equation}
From Lemma \ref{l:2}, we have that, for some positive constant $C_1$,
 \begin{align*}
 \frac1{C_1}|x|^{-(N-2-a_{\lambda_{\infty}})}&\leq\varphi_2(x)
\leq C_1|x|^{-(N-2-a_{\lambda_{\infty}})}\text{ in }\R^N\setminus B(0,2R_0).  
 \end{align*} 
Set $\Phi=\varphi_1+\varphi_2$. We claim that for $\tilde R>0$ large
enough there holds 
$$
-\Delta\Phi(x)-V(x)\Phi(x)-\frac{\gamma_{\infty}\alchi_{\R^N\setminus
    B(0,\tilde R)}(x)}{|x|^2}\Phi(x)>0,\quad\text{ a.e. in }\R^N. 
$$
For all $\tilde R>2 R_0$, from (\ref{eq:26}--\ref{eq:27}) we deduce 
\begin{align}
\label{eq:29}&-\Delta\Phi(x)-V(x)\Phi(x)-\frac{\gamma_{\infty}\alchi_{\R^N\setminus
    B(0,\tilde R)}(x)}{|x|^2}\Phi(x)\\[5pt] 
&\notag \geq\bigg[\Big(\frac{N-2}{2}\Big)^2\!\!-\e^2-\gamma_{\infty}\bigg]\frac{\alchi_{\R^N\setminus
    B(0,\tilde R)}(x)}{|x|^2}\varphi_1(x)+\bigg[\Big(\frac{N-2}{2}\Big)^2-\e^2\bigg]\frac{\alchi_{B(0,\tilde
    R)\setminus B(0,2R_0)}(x)}{|x|^2}\varphi_1(x)\\[5pt]  
&\notag \quad+f_1(x)+f_2(x)+
\bigg[W(x)\alchi_{B(0,2R_0)}(x)+\frac{\lambda_{\infty}\alchi_{B(0,2
  R_0)\setminus B(0, R_0)}(x)}{|x|^2}\bigg]\varphi_1(x)\\[5pt]
&\notag \quad-\frac{\gamma_{\infty}\alchi_{\R^N\setminus B(0,\tilde
    R)}(x)}{|x|^2}\varphi_2(x)-V(x)\varphi_1(x)\\[5pt] 
&\notag \geq f_1(x)+f_2(x)+\bigg[\Big(\frac{N-2}{2}\Big)^2\!\!-\e^2-\gamma_{\infty}-\lambda_{\infty}\bigg]\frac{\alchi_{\R^N\setminus
    B(0,\tilde R)}(x)}{|x|^2}\frac{1}{|x|^{\frac{N-2}{2}+\e}}\\[5pt] 
&\notag \quad+\bigg[\Big(\frac{N-2}{2}\Big)^2-\e^2-\lambda_{\infty}
-\frac{C_0}{|x|^{\delta}}\bigg]\frac{\alchi_{B(0,\tilde
R)\setminus B(0,2R_0)}(x)}{|x|^2}\varphi_1(x)\\  
&\notag \
-\frac{C_0}{|x|^{2+\delta}}\frac{1}{|x|^{\frac{N-2}{2}+\e}}\alchi_{\R^N\setminus
  B(0,\tilde R)}(x)
-\frac{\gamma_{\infty}}{|x|^2}\frac{C_1}{|x|^{N-2-a_{\lambda_\infty}}}\alchi_{\R^N\setminus
  B(0,\tilde R)}(x). 
\end{align}
By the choice of $f_2$, we have that in $B(0,2R_0)$ 
$$
-\Delta\Phi(x)-V(x)\Phi(x)-\frac{\gamma_{\infty}\alchi_{\R^N\setminus
    B(0,\tilde R)}(x)}{|x|^2}\Phi(x)\geq f_1(x)+f_2(x)>0.
$$
From (\ref{eq:29}) and the choice of $\e$ and $R_0$, it follows that, in
$B(0,\tilde R)\setminus B(0,2R_0)$,
\begin{align*}
-\Delta\Phi(x)-V(x)\Phi(x)-\frac{\gamma_{\infty}\alchi_{\R^N\setminus
    B(0,\tilde R)}(x)}{|x|^2}\Phi(x)\geq\bigg[\Big(\frac{N-2}{2}\Big)^2-\e^2-\lambda_{\infty}
-\frac{C_0}{|x|^{\delta}}\bigg]\frac{1}{|x|^2}\varphi_1(x)
>0.
\end{align*}
From (\ref{eq:29}) and the choice of $\e$, we deduce that, in $\R^N\setminus B(0,\tilde R)$,
\begin{align*}
-\Delta\Phi(x)&-V(x)\Phi(x)-\frac{\gamma_{\infty}\alchi_{\R^N\setminus
    B(0,\tilde R)}(x)}{|x|^2}\Phi(x)\geq\\
&\notag \geq \bigg[\Big(\frac{N-2}{2}\Big)^2\!\!-\e^2-\gamma_{\infty}-\lambda_{\infty}\bigg]\frac{1}{|x|^{2+\frac{N-2}{2}+\e}}
-\frac{C_0}{|x|^{2+\delta+\frac{N-2}{2}+\e}}
-\frac{\gamma_{\infty}C_1}{|x|^{N-a_{\lambda_\infty}}}>0
\end{align*}
provided $\tilde R$ is large enough. The claim is thereby proved.
\end{pf}

\begin{Theorem}\label{t:pertinf}
Let 
$$
V(x)=\sum_{i=1}^k\frac{\lambda_i \alchi_{B(a_i,r_i)}(x)}{|x-a_i|^2}+\frac{\lambda_{\infty} \alchi_{\R^N\setminus B(0,R)}(x)}{|x|^2}+W(x)\in{\mathcal V},
$$
where $W\in L^{\infty}(\R^N)$, $W(x)=O(|x|^{-2-\delta})$, with $\delta>0$, as $|x|\to\infty$. Assume that $\mu(V)>0$, namely that the quadratic form associated to the operator $-\Delta-V$ is positive definite.  Let $\gamma_{\infty}\in\R$ such that
$\lambda_{\infty}+\gamma_{\infty}<\big(\frac{N-2}2\big)^2$. Then there exists 
$\tilde R>R$ such that $\mu\big(V+\frac{\gamma_{\infty}}{|x|^2}\alchi_{\R^N\setminus B(0,\tilde R)}\big)>0$, namely the quadratic form associated to the operator $-\Delta-V-\frac{\gamma_{\infty}}{|x|^2}\alchi_{\R^N\setminus B(0,\tilde R)}$ is positive definite.
\end{Theorem}

\begin{pf}
As already observed in the proof of Lemma \ref{l:positivity_condition}, if $V\in{\mathcal V}$ and $\mu(V)>0$, then there exists $\e>0$ such that $V+\e V^+\in{\mathcal V}$, $\mu(V+\e V^+)>0$, and
$\lambda_{\infty}+\gamma_{\infty}+\e(\lambda_{\infty}^++\gamma_{\infty}^+)<\big(\frac{N-2}2\big)^2$. In particular $V+\e V^+$ also satisfies the assumptions of Lemma \ref{l:pertinf}, hence there exist $\tilde R>R$ and $\Phi\in\Di$, $\Phi>0$  in
$\R^N\setminus\{a_1,\dots,a_k\}$ and smooth in $\R^N\setminus\{a_1,\dots,a_k\}$ such that 
$$
-\Delta\Phi(x)-V(x)\Phi(x)-\frac{\gamma_{\infty}}{|x|^2}\alchi_{\R^N\setminus
  B(0,\tilde
  R)}\Phi(x)>\e\Big(V^+(x)+\frac{\gamma_{\infty}^+}{|x|^2}\Big)\Phi(x)\geq 
\e\Big(V(x)+\frac{\gamma_{\infty}}{|x|^2}\Big)^+\Phi(x).
$$
The conclusion follows now from Lemma \ref{l:positivity_condition}.
\end{pf}

\section{Separation Theorem}\label{sec:scattering}

In this section we provide a tool to construct a positive operator
from two positive potentials in $\mathcal V$ whose interaction at
infinity is not too strong. 
To this aim we first show how, starting from the supersolutions
corresponding to each positive given operator, it is possible to
scatter the singularities and obtain a positive supersolution for the
resulting operator by summation. 

\begin{Lemma}\label{l:scattering}
Let 
\begin{align*}
&V_1(x)=\sum_{i=1}^{k_1}\frac{\lambda_i^1\alchi_{B(a_i^1,r_i^1)}(x)}{|x-a_i^1|^2}+
\frac{\lambda_{\infty}^1\alchi_{\R^N\setminus
    B(0,R_1)}(x)}{|x|^2}+W_1(x)\in{\mathcal V},\\ 
&V_2(x)=\sum_{i=1}^{k_2}\frac{\lambda_i^2\alchi_{B(a_i^2,r_i^2)}(x)}{|x-a_i^2|^2}+
\frac{\lambda_{\infty}^2\alchi_{\R^N\setminus
    B(0,R_2)}(x)}{|x|^2}+W_2(x)\in{\mathcal V},
\end{align*}
where $W_i\in L^{\infty}(\R^N)$, $W_i(x)=O(|x|^{-2-\delta})$, $i=1,2$, with
$\delta>0$, as $|x|\to\infty$. Assume that $\mu(V_1),\mu(V_2)>0$, namely that the quadratic forms
associated to the operators $-\Delta-V_1$, $-\Delta-V_2$ are positive
definite and that $\lambda_{\infty}^1+\lambda_{\infty}^2<(\frac{N-2}{2})^2$. 
Then, there exists $R>0$ such that, for every $y\in\R^N$ with $|y|\geq
R$, there exists $\Phi_y\in\Di$, $\Phi_y\geq 0$ a.e. in $\R^N$,  $\Phi_y> 0$ 
and $C^1$-smooth in
$\R^N\setminus\{a_i^1,a_i^2+y\}_{i=1,\dots k_j,j=1,2}$, such that 
$$
-\Delta\Phi_y(x)-\left(V_1(x)+V_2(x-y)\right)\Phi_y(x)>0\quad\text{ a.e. in }\R^N.
$$
\end{Lemma}
 
\begin{pf}
Let $0<\e<\!\!<1$ be such that $\lambda_{\infty}^1+\lambda_{\infty}^2<(\frac{N-2}{2})^2-\e$ and, for $j=1,2$,
set 
$$
\Lambda=\Big(\frac{N-2}{2}\Big)^2-\e\quad\text{ and }\quad\gamma_{\infty}^j=\Lambda-\lambda_{\infty}^j.
$$
Let us also choose $0<\eta<\!\!<1$ such that 
\begin{equation}\label{eq:13}
\lambda^2_{\infty}<\gamma^1_{\infty}(1-2\eta)\quad\text{and}\quad
\lambda^1_{\infty}<\gamma^2_{\infty}(1-2\eta).
\end{equation}
We can choose $\bar R>0$ such that, for $j=1,2$, $\bigcup_{i=1}^{k_j}B(a_i^j,r_i^j)\subset B(0,\bar R)$,
and define
\begin{equation}\label{eq:scatt9}
p_j(x):=\begin{cases}
|x-a_i^j|^{-2+\sigma}&\text{in }B(a_i^j,r_i^j),\\
1&\text{in }B(0,\bar R)\setminus\bigcup_{i=1}^{k_j}B(a_i^j,r_i^j),\\
0&\text{in }\R^N\setminus B(0,\bar R),
\end{cases}
\end{equation}
with $\sigma>0$. In view of Theorem \ref{t:pertinf}, there exist $\tilde R_j$ such that the quadratic forms associated to the operators $-\Delta-V_j-\frac{\gamma_{\infty}^j}{|x|^2}\alchi_{\R^N\setminus
B(0,\tilde R_j)}$ are positive definite. 
Therefore, since $p_j\in L^{N/2}$, the infima 
$$
\mu_j=\inf_{u\in\Di\setminus\{0\}}\frac{\int_{\R^N}\big[|\n u(x)|^2-V_j(x)u^2(x)-\gamma_{\infty}^j|x|^{-2}\alchi_{\R^N\setminus
B(0,\tilde R_j)}u^2(x)\big]\,dx}{\int_{\R^N}p_j(x)u^2(x)\,dx}
$$
are achieved by some $\psi_j\in\Di$, $\psi_j\geq0$ a.e. in $\R^N$, $\psi_j>0$ in
$\R^N\setminus\{a_1^j,\dots,a_{k_j}^j\}$, solving equation
\begin{equation}\label{eq:8}
-\Delta\psi_j(x)-V_j(x)\psi_j(x)= \mu_{j}p_j(x)\psi_j(x)+\frac{\gamma_{\infty}^j}{|x|^2}\alchi_{\R^N\setminus
B(0,\tilde R_j)}\psi_j(x)\quad\text{ in }\R^N.
\end{equation}
In view of Lemma \ref{l:2}, there holds
$$
\lim_{|x|\to+\infty} \psi_j(x)|x|^{N-2-a_{\Lambda}}=\ell_j>0,
$$
hence the function $\varphi_j:=\frac{\psi_j }{\ell_j}$ solves (\ref{eq:8})
and $\varphi_j(x)\sim|x|^{-\left(N-2-a_{\Lambda}\right)}$ at $\infty$. 
Then there exists $\rho>\max\{\tilde R_1,\tilde R_2,\bar R\}$ such that, in  $\R^N\setminus B(0,\rho)$, 
\begin{equation}\label{eq:14}
(1-\eta^2)|x|^{-\left(N-2-a_{\Lambda}\right)}\leq\varphi_j(x)\leq
(1+\eta)|x|^{-\left(N-2-a_{\Lambda}\right)} 
\end{equation}
and that 
\begin{equation}\label{eq:15}
|W_1(x)|\leq\eta\gamma_\infty^2|x|^{-2}\quad\mbox{ and }\quad|W_2(x)|\leq\eta\gamma_\infty^1|x|^{-2}.
\end{equation}
Moreover, from Lemma \ref{l:1} we can deduce 
that for some positive constant $C$
\begin{align}\label{eq:16}
\frac1C |x-a_i^j|^{-a_{\lambda_i^j}}&\leq\varphi_j(x)\leq C|x-a_i^j|^{-a_{\lambda_i^j}}\text{ in }B(a_i^j,r_i^j),\quad i=1,\dots,k_j
\end{align}
and $\varphi_j$ are $C^1$-smooth outside the poles. 

For any $y\in\R^N$, let us consider the function
$$
\Phi_y(x):=\gamma_{\infty}^2\varphi_1(x)+\gamma_{\infty}^1\varphi_2(x-y)\in\Di.
$$
Then 
\begin{align*}
-\D\Phi_y(x)-&\big(V_1(x)+V_2(x-y)\big)\Phi_y(x)\\
&=\mu_{1}\gamma_{\infty}^2 p_1(x)\varphi_1(x)+
\frac{\gamma_{\infty}^1\gamma_{\infty}^2}{|x|^2}\alchi_{\R^N\setminus
B(0,\tilde R_1)}\varphi_1(x)\\
&\quad+\mu_{2}\gamma_{\infty}^1 p_2(x-y)\varphi_2(x-y)+ \frac{\gamma_{\infty}^1\gamma_{\infty}^2}{|x-y|^2}\alchi_{\R^N\setminus B(y,\tilde R_2)}\varphi_2(x-y)\\
&\quad-\gamma_{\infty}^1V_1(x)\varphi_2(x-y)- \gamma_{\infty}^2V_2(x-y)\varphi_1(x).
\end{align*}
From (\ref{eq:13}) and (\ref{eq:14}), it follows that in $\R^N\setminus\big(B(0,\rho)\cup B(y,\rho)\big)$  
\begin{align*}
-\D\Phi_y(x)-&\big(V_1(x)+V_2(x-y)\big)\Phi_y(x)\\
&\geq\frac{\gamma_{\infty}^1\gamma_{\infty}^2}{|x|^2}\varphi_1(x)+
\frac{\gamma_{\infty}^1\gamma_{\infty}^2}{|x-y|^2}\varphi_2(x-y)\\ 
&\quad-\gamma_{\infty}^1\bigg(\frac{\lambda_{\infty}^1}{|x|^2}+W_1(x)\bigg)
\varphi_2(x-y)-\gamma_{\infty}^2\bigg(\frac{\lambda_{\infty}^2}{|x-y|^2}+
W_2(x-y)\bigg)\varphi_1(x)\\
&>\frac{\gamma_{\infty}^1\gamma_{\infty}^2(1-\eta^2)}{|x|^{N-a_{\Lambda}}}+
\frac{\gamma_{\infty}^1\gamma_{\infty}^2(1-\eta^2)}{|x-y|^{N-a_{\Lambda}}}-
\frac{\gamma_{\infty}^1\gamma_{\infty}^2(1-\eta^2)}{|x|^2|x-y|^{N-2-a_{\Lambda}}}-
\frac{\gamma_{\infty}^1\gamma_{\infty}^2(1-\eta^2)}{|x-y|^2|x|^{N-2-a_{\Lambda}}}\\
&=\gamma_{\infty}^1\gamma_{\infty}^2(1-\eta^2)\left(\frac{1}{|x|^{N-a_{\Lambda}-2}}-
  \frac{1}{|x-y|^{N-a_{\Lambda}-2}}\right)\left(\frac{1}{|x|^2}-
  \frac{1}{|x-y|^2}\right)\geq0. 
\end{align*}

For $|y|$ sufficiently large, $B(0,\rho)\cap B(y,\rho)=\emptyset$. In $B(a_i^1,r_i^1)$, for $i=1,\dots,k_1$, from (\ref{eq:scatt9}), (\ref{eq:14}), (\ref{eq:15}), (\ref{eq:16}) we have that 
\begin{align*}
-\D\Phi_y(x)-&\big(V_1(x)+V_2(x-y)\big)\Phi_y(x)\\
&\geq\mu_{1}\gamma_{\infty}^2p_1(x)\varphi_1(x)+ \frac{\gamma_{\infty}^1\gamma_{\infty}^2}{|x-y|^2}\varphi_2(x-y)\\
&\quad-\gamma_{\infty}^1V_1(x)\varphi_2(x-y)- \gamma_{\infty}^2V_2(x-y)\varphi_1(x)\\
&\geq|x-a_i^1|^{-a_{\lambda_i^1}-2+\sigma}\left[\frac{\mu_{1}\gamma_{\infty}^2}C+o(1)\right], \quad\mbox{ as } |y|\to\infty.
\end{align*}
In $B(0,\rho)\setminus\bigcup_{i=1}^{k_1}B(a_i^1,r_i^1)$, from (\ref{eq:scatt9}), (\ref{eq:14}), (\ref{eq:15}) and since $\varphi_1>c>0$, we obtain that
$$
-\D\Phi_y(x)-\big(V_1(x)+V_2(x-y)\big)\Phi_y(x)\geq \mu_{1}\gamma_{\infty}^2c+o(1),\quad\mbox{ as }|y|\to\infty.
$$

In a similar way we can prove that, if $|y|$ is large enough, 
$$
-\D\Phi_y(x)-\big(V_1(x)+V_2(x-y)\big)\Phi_y(x)\geq0,\quad\mbox{ a.e. in }B(y,\rho).
$$
\end{pf}

\begin{pfn}{Theorem \ref{t:scattering}}
Let us fix $\e\in(0,1)$ such that, for $j=1,2$, 
\begin{align}
\e<\min\bigg\{
2S\,\mu(V_j),  \frac{\mu(V_j)}{4}\bigg[\frac4{(N-2)^2}
\Big(\sum_{i=1}^k(\lambda_i^j)^++(\lambda^j_{\infty})^+\Big)
+S^{-1}\|W_j\|_{L^{N/2}(\R^N)}\bigg]^{-1}\bigg\}\label{eq:18}
\end{align}
and
\begin{align*}
\mu(V_j+\e V_j^+)>0,\quad\lambda_{\infty}^1+\lambda_{\infty}^2+\e\left(\left(\lambda_{\infty}^1\right)^++
  \left(\lambda_{\infty}^2\right)^+\right)<\Big(\frac{N-2}2\Big)^2,
\end{align*}
see the proof of Lemma \ref{l:positivity_condition}. Fix $R>0$ such
that
\begin{align}\label{eq:23}
\|W_j\alchi_{\R^N\setminus
  B(0,R)}\|_{L^{N/2}(\R^N)}<\min\Big\{\frac\e4,\frac\e{16}\,S\Big\}.
\end{align}
Denoting $V_{j,R}:=V_j-W_j\alchi_{\R^N\setminus B(0,R)}$, from
(\ref{eq:18}) and (\ref{eq:23}), there results
\begin{multline*}
\int_{\R^N}\big(|\n u(x)|^2-V_{j,R}(x)u^2(x)\big)\,dx\\
\geq \Big[\mu(V_j)-\|W_j\alchi_{\R^N\setminus
  B(0,R)}\|_{L^{N/2}(\R^N)}S^{-1}\Big]
\int_{\R^N}|\n u(x)|^2\,dx\geq \frac{\mu(V_j)}2\int_{\R^N}|\n
u(x)|^2\,dx,
\end{multline*}
therefore, from (\ref{eq:18}), it follows
\begin{align*}
\int_{\R^N}&\big(|\n u(x)|^2-\big(V_{j,R}(x)+\e V_{j,R}^+(x)\big)
u^2(x)\big)\,dx\\
&\geq \left[\frac{\mu(V_j)}2-\e \bigg(\frac4{(N-2)^2}
\Big(\sum_{i=1}^k(\lambda_i^j)^++(\lambda^j_{\infty})^+\Big)
+S^{-1}\|W_j\|_{L^{N/2}(\R^N)}\bigg)\right]\int_{\R^N}|\n
u(x)|^2\,dx\\
&\geq \frac{\mu(V_j)}4\int_{\R^N}|\n
u(x)|^2\,dx.
\end{align*}
Hence the potentials $V_{j,R}+\e V_{j,R}^+$ satisfy the assumptions of
Lemma \ref{l:scattering}, which yields, for $|y|$ sufficiently large,
the existence of $\Phi_y\in\Di$, $\Phi_y\geq 0$ a.e. in $\R^N$,  $\Phi_y> 0$
and $C^1$-smooth in
$\R^N\setminus\{a_i^1,a_i^2+y\}_{i=1,\dots k_j,j=1,2}$, such that 
\begin{equation*}
-\Delta\Phi_y(x)-V_{R,y}(x)\Phi_y(x)>
\e\, V_{R,y}^+(x)\Phi_y(x)
\quad\text{ a.e. in }\R^N,
\end{equation*}
where $V_{R,y}(x):=V_{1,R}(x)+V_{2,R}(x-y)$. As a consequence, arguing as in the proof
of Lemma \ref{l:positivity_condition}, we easily deduce that 
\begin{equation}\label{eq:25}
\inf_{u\in\Di\setminus\{0\}}\frac{\int_{\R^N}\big(|\n
u(x)|^2-V_{R,y}(x)u^2(x)\big)\,dx}{\int_{\R^N}V_{R,y}^+(x)u^2(x)\,dx}\geq\e.
\end{equation}
We claim that 
$$
\inf_{u\in\Di\setminus\{0\}}\frac{\int_{\R^N}\big(|\n
u(x)|^2-\big(V_1(x)+V_2(x-y)\big)
u^2(x)\big)\,dx}{\int_{\R^N}|\n u(x)|^2\,dx}>0.
$$
We argue by contradiction and assume that
there exists a sequence $u_n\in\Di$ such that 
$$
\lim_{n\to\infty}\int_{\R^N}\big(|\n u_n(x)|^2-\big(V_1(x)+V_2(x-y)\big)
u^2_n(x)\big)\,dx=0\quad\text{ and }\quad 
\int_{\R^N}|\n u_n(x)|^2=1.
$$
From (\ref{eq:25}) and (\ref{eq:23}), we obtain, for $n$ large enough,
\begin{align*}
\e&=
\e\int_{\R^N}\!\!\Big[  V^+_{R,y}(x)+  W_1(x)\alchi_{\R^N\setminus B(0,R)}(x)
+  W_2(x-y)\alchi_{\R^N\setminus B(y,R)}(x)- 
V^-_{R,y}(x)\Big]u_n^2(x)\,dx +o(1)\\
&\leq \int_{\R^N}|\n u_n(x)|^2\,dx-\int_{\R^N}
V_{R,y}(x)u_n^2(x)\,dx+\e S^{-1}\sum_{j=1}^2
\|W_j\alchi_{\R^N\setminus
  B(0,R)}\|_{L^{N/2}(\R^N)} +o(1)\\
&\leq \int_{\R^N}\big(|\n
u_n(x)|^2-\big(V_1(x)+V_2(x-y)\big)
u_n^2(x)\big)\,dx+{\textstyle{\frac{(1+\e)}S}}\sum_{j=1}^2
\|W_j\alchi_{\R^N\setminus
  B(0,R)}\|_{L^{N/2}(\R^N)} +o(1)\\
&\leq \frac\e8(1+\e)+o(1)\leq \frac\e2,
\end{align*}
which is a contradiction. The theorem is thereby proved.
\end{pfn}

\section{Positivity of multipolar operators}\label{mainresult}

This section is devoted to the proof of Theorem \ref{t:mainresult}. 

\bigskip\noindent

\begin{pfn}{Theorem \ref{t:mainresult}}
\mbox{}

\medskip\noindent
{\bf Step 1 (Sufficiency.)}\quad
We prove the sufficient condition on the positivity of the quadratic
form $Q_{\lambda_1,\dots,\lambda_k,a_1,\dots,a_k}$ applying an iterating process
on the number of poles $k$. Let $\lambda_1\leq\lambda_2\leq\cdots\leq\lambda_k$.

As observed in Section \ref{intro}, if $k=2$ the claim is true for any choice of $a_1,a_2$.
Suppose that the claim is true for $k-1$, let us prove it for $k$. 
 We may assume $\lambda_k>0$,
otherwise the proof is trivial. If $\lambda_1,\dots,\lambda_k$ satisfy
(\ref{eq:19}), then the same holds true for $\lambda_1,\dots,\lambda_{k-1}$. By
the recursive assumption, there exists a configuration of poles
$\{a_1,\dots,a_{k-1}\}$ such that the quadratic form
$Q_{\lambda_1,\dots,\lambda_{k-1},a_1,\dots,a_{k-1}}$ associated to the operator

\begin{equation*}
L_{\lambda_1,\dots,\lambda_{k-1},a_1,\dots,a_{k-1}}=
-\D-{\displaystyle{\sum_{i=1}^{k-1}}}\dfrac{\lambda_i}{|x-a_i|^2} 
\end{equation*}
is positive definite.  

We claim that there exists $a_k\in\R^N$ such that the quadratic form
associated to the operator $L_{\lambda_1,\dots,\lambda_k,a_1,\dots,a_k}$ is
positive definite. Indeed, the two potentials
$$
V_1(x)=\sum_{i=1}^{k-1}\frac{\lambda_i}{|x-a_i|^2},\quad V_2(x)=\frac{\lambda_k}{|x|^2},
$$
belong to the class $\mathcal V$ and satisfy the assumptions of
Theorem \ref{t:scattering}, which ensures the existence of $a_k\in\R^N$
such that the quadratic form associated to the operator
$$
L_{\lambda_1,\dots,\lambda_k,a_1,\dots,a_k}=-\Delta-\left(V_1+V_2(\cdot-a_k)\right)
$$
is positive definite.
    
\medskip\noindent
{\bf Step 2 (Necessity.)}\quad   
Assume that for some configuration $\{a_1,\dots, a_k\}$ and for some $\e>0$ there holds
$$
\int_{\R^N}|\n u(x)|^2dx-\sum_{i=1}^k {\lambda_i}\int_{\R^N}\frac{u^2(x)}{|x-a_i|^2}\,dx
\geq\e\int_{\R^N}|\n u(x)|^2\,dx,\quad\text{ for all }u\in\Di.
$$
Arguing by contradiction, suppose that, for some $i$, $\lambda_i\geq\big(\frac{N-2}2\big)^2$. Let $\delta\in(0,\e(N-2)^2/4)$. By
optimality of the best constant in the Hardy inequality (\ref{eq:hardy}) and by density of $C^{\infty}_{\rm c}(\R^N)$ in $\Di$, there exists $\phi\in C^{\infty}_{\rm c}(\R^N)$ such that 
$$
\int_{\R^N}|\n\phi(x)|^2dx-{\lambda_i}\int_{\R^N}\frac{\phi^2(x)}{|x|^2}\,dx
<\delta\int_{\R^N}\frac{\phi^2(x)}{|x|^2}\,dx.
$$
The rescaled function $\phi_{\mu}(x)=\mu^{-(N-2)/2}\phi(x/\mu)$ satisfies
\begin{align*}
Q_{\lambda_1,\dots,\lambda_k,a_1,\dots,a_k}(\phi_{\mu}(x-a_i))
&=\int_{\R^N}|\n\phi(x)|^2dx-{\lambda_i}\int_{\R^N}\frac{\phi^2(x)}{|x|^2}\,dx -\sum_{j\neq i}\lambda_j\int_{\R^N}\frac{\phi^2(x)}{\big|x-\frac{a_j-a_i}{\mu}\big|^2}\,dx\\
&=\int_{\R^N}|\n\phi(x)|^2dx- {\lambda_i}\int_{\R^N}\frac{\phi^2(x)}{|x|^2}\,dx+o(1),\quad\text{ as }\mu\to0.
\end{align*}
Letting $\mu\to0$, by Hardy's inequality, we obtain 
\begin{align*}
\e\int_{\R^N}|\n\phi(x)|^2dx &\leq\int_{\R^N}|\n\phi(x)|^2dx-{\lambda_i}\int_{\R^N}\frac{\phi^2(x)}{|x|^2}\,dx\\
&<\delta\int_{\R^N}\frac{\phi^2(x)}{|x|^2}\,dx
\leq\frac{4\delta}{(N-2)^2}\int_{\R^N}|\n \phi(x)|^2dx
\end{align*}
thus giving rise to a contradiction.

Suppose now that, $\Lambda:=\sum_{i=1}^k\lambda_i\geq\big(\frac{N-2}2\big)^2$. 
Let $\delta\in(0,\e(N-2)^2/4)$. As above, there exists $\phi\in
C^{\infty}_{\rm c}(\R^N)$ such that  
$$
\int_{\R^N}|\n\phi(x)|^2dx-{\Lambda}\int_{\R^N}\frac{\phi^2(x)}{|x|^2}\,dx
<\delta\int_{\R^N}\frac{\phi^2(x)}{|x|^2}\,dx.
$$
The rescaled function $\phi_{\mu}(x)=\mu^{-(N-2)/2}\phi(x/\mu)$ satisfies
\begin{align*}
Q_{\lambda_1,\dots,\lambda_k,a_1,\dots,a_k}(\phi_{\mu}(x))
&=\int_{\R^N}|\n\phi(x)|^2dx-
{\sum_{i=1}^k\lambda_i}\int_{\R^N}\frac{\phi^2(x)}{|x-a_i/\mu|^2}\,dx\\ 
&=\int_{\R^N}|\n\phi(x)|^2dx-
{\Lambda}\int_{\R^N}\frac{\phi^2(x)}{|x|^2}\,dx+o(1),\quad\text{ as
}\mu\to\infty, 
\end{align*} 
see \cite[Proposition 3.1]{FT}. Letting $\mu\to\infty$ and arguing as
above, we obtain easily a contradiction.~\end{pfn}

\section{Best constants in Hardy multipolar
inequalities} \label{sec:best-constants-hardy}

The classical Hardy's inequality states that 
$
\mu\big(|x|^{-2}\big)=1-\frac4{(N-2)^2}
$
is not attained. On the other hand, when dealing with multipolar
Hardy-type potentials, a balance between positive and negative
interactions between the poles can lead to attainability of the
best constant in the Hardy-type inequality associated to the
multisingular potential $V\in{\mathcal V}$:
\begin{equation}\label{eq:hardymul}
\int_{\R^N}V(x)\,|u(x)|^2\,dx\leq \big(1-\mu(V)\big)
\int_{\R^N}|\nabla u(x)|^2\,dx,\quad\text{ for all }u\in\Di.
\end{equation}
We recall that, in view of Lemma \ref{l:estmu}, the best constant in
inequality \eqref{eq:hardymul} can be estimated by 
terms of the best constant in the inequality associated to the
potential with one singularity located at the pole carrying
the largest mass, i.e. $\mu(V)\leq
1-\frac4{(N-2)^2}\,\max\{0,\lambda_1,\dots,\lambda_k,\lambda_{\infty}\}$.  We now prove
attainability of $\mu(V)$ when it 
stays strictly below the bound provided in Lemma~\ref{l:estmu}. 

\bigskip\noindent

\begin{pfn}{Proposition \ref{p:attai}}
Let us denote
$$
\bar\lambda=\max\{0,\lambda_1,\dots,\lambda_k,\lambda_{\infty}\}.
$$
From assumption  (\ref{eq:77}), there exists $\alpha>0$ such that
$\mu(V)=1-\frac{4\bar\lambda}{(N-2)^2}-\alpha$.
From Lemma \ref{l:sep}, there exists
$\d>0$ such that 
\begin{equation}\label{eq:76}
\mu(\widetilde V)\geq 1-\frac{4\bar\lambda}{(N-2)^2}-\frac{\alpha}2,
\end{equation}
where $\widetilde V(x)=
\sum_{i=1}^k \lambda_i\alchi_{B(a_i,\d)}|x-a_i|^{-2}+
\lambda_{\infty}\alchi_{\R^N\setminus B(0,R_0)}|x|^{-2}$. 
We notice that we can split $V=\widetilde V+\widetilde W$ for some
$\widetilde W\in L^{N/2}(\R^N)$. Let $u_n\in\Di$ be a minimizing sequence
for $\mu(V)$, namely  
$$
\int_{\R^N}|\n u_n(x)|^2\,dx=1\quad\text{and}\quad 
\int_{\R^N}\big(|\n u_n(x)|^2-V(x)u_n^2(x)\big)\,dx=\mu(V)+o(1)
\quad\text{as }n\to\infty.
$$
Being $\{u_n\}_n$ bounded in $\Di$, we can assume that, up to a
subsequence still denoted as $u_n$, $u_n$ converges to
some $u$ a.e. and weakly in $\Di$.   Since
$$
\mu(\widetilde V)\leq \int_{\R^N}\big(|\n
u_n(x)|^2-V(x)u_n^2(x)\big)\,dx+\int_{\R^N}\widetilde
W(x)u_n^2(x)\,dx
=\mu(V)+\int_{\R^N}\widetilde W(x)u^2(x)\,dx+o(1)
$$
as $n\to\infty$, from (\ref{eq:76}) and the choice of $\alpha$ it follows
that
$$
1-\frac{4\bar\lambda}{(N-2)^2}-\frac{\alpha}2\leq\mu(V)+\int_{\R^N}\widetilde W(x)
u^2(x)\,dx=1-\frac{4\bar\lambda}{(N-2)^2}-\alpha+\int_{\R^N}\widetilde W(x)
u^2(x)\,dx,
$$
hence $\int_{\R^N}\widetilde W(x)u^2(x)\,dx\geq\frac{\alpha}2>0$, thus implying
$u\not\equiv 0$. From weak convergence of $u_n$ to $u$, we 
deduce that
\begin{align*}
&\frac{\int_{\R^N}\big(|\n u(x)|^2-V(x)u^2(x)\big)\,dx}{\int_{\R^N}|\n
  u(x)|^2\,dx}\\[5pt]
&\ =\frac{\big[\int_{\R^N}\big(|\n
  u_n(x)|^2-V(x)u_n^2(x)\big)\,dx\big]-\big[ \int_{\R^N}\big(|\n
  (u_n-u)(x)|^2-V(x)(u_n-u)^2(x)\big)\,dx   \big]+o(1)}{\int_{\R^N}|\n
  u_n(x)|^2\,dx-\int_{\R^N}|\n
  (u_n-u)(x)|^2\,dx+o(1)}\\[5pt]
&\ \leq \mu(V)\,\frac{1-\int_{\R^N}|\n
  (u_n-u)(x)|^2\,dx+o(1)}{1-\int_{\R^N}|\n (u_n-u)(x)|^2\,dx+o(1)}=
\mu(V)\bigg[1+\frac{o(1)}{\int_{\R^N}|\n u(x)|^2\,dx+o(1)}\bigg]
\quad\text{as }n\to\infty.
\end{align*}
Letting $n\to\infty$, we  obtain that $u$ attains the infimum defining
$\mu(V)$.
\end{pfn}




\noindent As a consequence of the attainability of $\mu(V)$, a result of
continuity follows. 

\begin{Lemma}\label{l:cneg}
Let
\begin{align*}
&V(x)=\sum_{i=1}^{k}\frac{\lambda_i\alchi_{B(a_i,r_i)}(x)}{|x-a_i|^2}+
\frac{\lambda_{\infty}\alchi_{\R^N\setminus
B(0,R_0)}(x)}{|x|^2}+W(x)\in{\mathcal V},\\ 
&V_n(x)=\sum_{i=1}^{k}\frac{\lambda_i\alchi_{B(a_i^n,r_i)}(x)}{|x-a_i^n|^2}+
\frac{\lambda_{\infty}\alchi_{\R^N\setminus
B(0,R_0)}(x)}{|x|^2}+W_n(x)\in{\mathcal V}
\end{align*}
be such that $a_i^n\to a_i$ as $n\to\infty$, for all $i=1,\dots,k$, and $W_n\to W$  
in $L^{N/2}(\R^N)$.
Then 
$$
\lim_{n\to\infty}\mu(V_n)=\mu(V).
$$
\end{Lemma}

\begin{pf}
Let $\bar\lambda=\max\{0,\lambda_1,\dots,\lambda_k,\lambda_{\infty}\}$.
For any $u\in\Di$, Lemma \ref{l:cf} 
and the strong $L^{N/2}$-convergence of $W_n$ 
imply that 
$\int_{\R^N} V_n(x)u^2(x)\,dx\to\int_{\R^N}
V(x)u^2(x)\,dx$, hence
\begin{align*}
\mu(V_n)&\leq \frac{\int_{\R^N}\big(|\n
u(x)|^2- V_n(x)u^2(x)\big)\,dx}{\int_{\R^N}|\n u(x)|^2\,dx}=
\frac{\int_{\R^N}\big(|\n
u(x)|^2-V(x)u^2(x)\big)\,dx+o(1)}{\int_{\R^N}|\n u(x)|^2\,dx}
\end{align*}
for any $u\in \Di\setminus\{0\}$.
Therefore, letting $n\to\infty$ and taking infimum over $\Di\setminus\{0\}$, we
obtain that
\begin{equation}\label{eq:75}
\limsup_{n\to\infty}\mu(V_n)\leq\mu(V).
\end{equation}
In particular, the sequence $\big\{\mu(V_n)\big\}_n$ is bounded. We
now claim that 
\begin{equation}\label{eq:79}
\mu(V)= \liminf_{n\to\infty}\mu(V_n).
\end{equation}
Indeed, let $\big\{\mu(V_{n_j})\big\}_j$ be a subsequence such that
$\lim_j \mu(V_{n_j})=\liminf_{n\to\infty}\mu(V_n)$ and suppose, by
contradiction, that $\lim_j\mu(V_{n_j})<\mu(V)-\a$, for some $\a>0$.

In view of Lemma \ref{l:sep}, see also Remark
\ref{r:shattering-lemma}, there exists $\d>0$ independent of $n$ such
that
\begin{equation}\label{eq:ne}
\mu(\widetilde
V_{n})\ \geq 1-\frac{4\bar\lambda}{(N-2)^2}-\frac{\alpha}2,\quad
\mu(\widetilde
V) \geq 1-\frac{4\bar\lambda}{(N-2)^2}-\frac{\alpha}2,
\end{equation}
where
\begin{align*}
\widetilde
V_n(x)=\sum_{i=1}^k\frac{\lambda_i\alchi_{B(a_i^n,\d)}}{|x-a_i^n|^2}+
\frac{\lambda_{\infty}\alchi_{\R^N\setminus B(0,R_0)}}{|x|^2},
\quad\widetilde V(x)=
\sum_{i=1}^k\frac{\lambda_i\alchi_{B(a_i,\d)}}{|x-a_i|^2}+
\frac{\lambda_{\infty}\alchi_{\R^N\setminus B(0,R_0)}}{|x|^2}. 
\end{align*}
We can write $V_n=\widetilde V_n+\widetilde W_n$ and $V=\widetilde
V+\widetilde W$, where $$\widetilde
W_n=W_n+\sum_{i=1}^k\frac{\lambda_i}{|x-a_i^n|^2}\alchi_{B(a_i^n,r_i)\setminus
  B(a_i^n,\d)}\quad\text{ and }\quad\widetilde
W=W+\sum_{i=1}^k\frac{\lambda_i}{|x-a_i|^2}\alchi_{B(a_i,r_i)\setminus
  B(a_i,\d)}.$$  
By the Dominated Convergence Theorem we deduce that $\widetilde
W_n\to\widetilde W$ in $L^{N/2}(\R^N)$.

From Lemma \ref{l:estmu}, we have that, for large $j$, 
\begin{equation}\label{eq:78}
\mu(V_{n_j})<\mu(V)-\a< 1-\frac{4\bar\lambda}{(N-2)^2},
\end{equation}
hence, by Proposition \ref{p:attai}, $\mu(V_{n_j})$ is attained by some
$\varphi_j\in\Di$ satisfying   
$$
\int_{\R^N}|\n \varphi_j(x)|^2\,dx=1\quad\text{and}\quad
\int_{\R^N}\big(|\n
\varphi_j(x)|^2-V_{n_j}(x)\varphi_j^2(x)\big)\,dx=\mu(V_{n_j}).
$$
Moreover $\varphi_j$ satisfies the equation
\begin{equation}\label{eq:81}
-\D\varphi_j(x)-V_{n_j}(x)\varphi_j(x)=-\mu(V_{n_j})\D\varphi_j(x).
\end{equation}
Since $\{\varphi_j\}_j$ is bounded in $\Di$, there exists a
subsequence, still denoted as $\{\varphi_j\}_j$, weakly converging to
some $\varphi$ in $\Di$.
From \eqref{eq:78} and \eqref{eq:ne} it follows that
\begin{align*}
1-\frac{4\bar\lambda}{(N-2)^2}-\frac{\alpha}2
&\leq\mu(\widetilde V_{n_j})\leq \int_{\R^N}\big(|\n
\varphi_j(x)|^2-V_{n_j}(x)\varphi_j^2(x)\big)\,dx+\int_{\R^N}\widetilde W_{n_j}(x)\varphi_j^2(x)\,dx\\
&=\mu(V_{n_j})+\int_{\R^N}\widetilde W_{n_j}(x)\varphi_j^2(x)\,dx
<\mu(V)-\a+\int_{\R^N}\widetilde W(x)\varphi^2(x)\,dx+o(1),
\end{align*}
as $j\to\infty$. Letting $j\to\infty$, we obtain that
\begin{gather*}
1-\frac{4\bar\lambda}{(N-2)^2}-\frac{\alpha}2\leq\mu(V)
-\a+\int_{\R^N}\!\widetilde W(x)\varphi^2(x)\,dx
\leq
1-\frac{4\bar\lambda}{(N-2)^2}-\alpha+\int_{\R^N}\!\widetilde W(x)\varphi^2(x)\,dx,
\end{gather*}
yielding  $\int_{\R^N}\widetilde W(x)\varphi^2(x)\,dx\geq\frac{\alpha}2>0$ and hence
$\varphi\not\equiv0$. 
We claim that 
\begin{equation}\label{eq:82}
\lim_{j\to\infty}\int_{\R^N}\widetilde V_{n_j}(x)\varphi_j(x)v(x)\,dx=
\int_{\R^N}\widetilde V(x)\varphi(x)v(x)\,dx\quad\text{for all
}v\in\Di.
\end{equation}
Indeed for any $v\in\Di$ and $\e>0$, by density there exists $\psi\in C^{\infty}_{\rm 
  c}(\R^N\setminus\{a_1,\dots,a_k\})$ such that
$\|v-\psi\|_{\Di}<\e$. Since $\psi$ lies far away from the
singularities, from Hardy's inequality we have that
\begin{multline*}
\bigg|\int_{\R^N}\widetilde V_{n_j}(x)\varphi_j(x)v(x)\,dx-
\int_{\R^N}\widetilde V(x)\varphi(x)v(x)\,dx\bigg|\\
\leq {\rm const\,}\e+
\bigg|
\int_{\R^N}\big(\widetilde V_{n_j}-\widetilde
V)(x)\varphi_j(x)\psi(x)\,dx\bigg|
\leq {\rm const\,}\e+o(1)\quad\text{as }j\to\infty.
\end{multline*}
(\ref{eq:82}) is thereby proved. From (\ref{eq:82}) and strong $L^{N/2}$-convergence of $\widetilde W_n$
to $\widetilde W$, we can multiply (\ref{eq:81}) by $\varphi$ and pass to
limit as $j\to\infty$ thus obtaining
$$
\int_{\R^N}\big(|\n
\varphi(x)|^2-V(x)\varphi^2(x)\big)\,dx=\liminf_{n\to\infty}\mu(V_{n})\int_{\R^N}|\n
\varphi(x)|^2\,dx,
$$
and consequently
$\liminf_{n\to\infty}\mu(V_n)\geq\mu(V)$, a contradiction.
Claim (\ref{eq:79}) is thereby proved. The conclusion follows from
(\ref{eq:75}) and (\ref{eq:79}).
\end{pf}

\begin{remark}\label{r:n2cont}
We emphasize that $\mu: L^{N/2,\infty}(\R^N)\to\R$, $\mu:\ V\mapsto
\mu(V)$ is continuous with respect to the $L^{N/2,\infty}$-norm. In
particular the first $\Di$-eigenvalue $\mu(V)$ is continuous not only
with respect to the location of the singularities but also with
respect to their masses $\lambda_i$'s.  
\end{remark}

\begin{remark}\label{r:nocont}
We notice that if $V_n\in{\mathcal V}$ converge to
$V\in{\mathcal V}$ in the sense that the poles of $V_n$ converge to
the poles of $V$ (i.e. in the sense of Lemma \ref{l:cneg}),
then $\|V_n-V\|_{L^{N/2,\infty}(\R^N)}$ does
not tend to zero. On the other hand $\mu(V_n)\to\mu(V)$. In other
words, the first $\Di$-eigenvalue is stable with respect to small
perturbations of configurations of poles, even though  such
perturbations make the $L^{N/2,\infty}$-distance between the
potentials far away from zero.  
\end{remark}

\begin{remark}\label{rem:segno}
  Since operators with potentials in ${\mathcal V}$ are compact
  perturbations of positive operators, see Lemma \ref{l:semibounded},
  we can prove, for the first weighted eigenvalue $\nu_p(V)$ defined
  in (\ref{eq:59}), the same type of continuity result as in Lemma
  \ref{l:cneg}. Moreover, it is not difficult to check that the
  infimum defining $\nu_p$ is attained, for any $p\in
  L^{N/2}(\R^N)\cap C^0(\R^N)$, $p>0$. As a consequence, for any
  $V\in{\mathcal V}$, there holds
$$
\sgn \nu_p(V)=\sgn \mu(V),\quad\text{where}\quad \sgn t:=
\begin{cases}
1,&\text{if }t>0,\\
0, &\text{if }t=0,\\
-1,&\text{if }t<0.
\end{cases}
$$
Indeed, it is obvious that $\nu_p(V)<0$ if and only if $\mu(V)<0$. If
$\mu(V)>0$, from H\"older's and Sobolev's inequalities it immediately
follows that $\nu_p(V)>0$. Assume now that $\nu_p(V)>0$. By
continuity, we can find $\e>0$ such that $\nu_p(V+\e V^+)>0$ is
attained, thus providing a positive function satisfying condition
$(ii)$ of Lemma \ref{l:positivity_condition}, which yields the
positivity of~$\mu(V)$.
\end{remark}

\section{Essential self-adjointness}\label{sec:semi-bound-essent}

The Shattering Lemma \ref{l:sep} reveals how
Schr\"odinger operators with potentials lying in the class ${\mathcal
  V}$ are actually compact perturbations of positive operators,
see Lemma \ref{l:semibounded}. Hence they are {\em semi-bounded
  symmetric operators} and their 
 $L^2(\R^N)$-spectrum  is bounded from below. Consequently 
the class ${\mathcal V}$ provides us with a good framework  
to study the spectral properties of multisingular  Schr\"odinger
operators in $L^2(\R^N)$.

For any $V\in\mathcal V$,
let us discuss essential self-adjointness of the operator $-\D-V$ on the domain
$C^{\infty}_{\rm c}(\R^N\setminus\{a_1,\dots,a_k\})$. In the case of
just one singularity (i.e. $k=1$), a complete answer to this problem
is contained in a theorem due to Kalf,  Schmincke, Walter, and W\"ust
\cite{kwss}  and Simon \cite{simon73}
(see also \cite[Theorems X.11 and X.30]{reedsimon}):
\begin{Theorem}{\bf [\,Kalf,  Schmincke, Walter,  W\"ust, Simon\,]}\label{t:kwss}
Let $V(x)=\frac{\lambda}{|x|^2}+W(x)$, $W\in L^{\infty}(\R^N)$.
The operator $-\D-V$ is essentially self-adjoint
in $C^{\infty}_{\rm c}(\R^N\setminus\{0\})$ if and only if $\lambda\leq
(N-2)^2/4-1$. 
\end{Theorem} 
We are now going to extend the above result to potentials lying in
the class ${\mathcal V}$, for which we give below a 
{\em{self-adjointness  criterion}}. According to Lemma \ref{l:semibounded}, 
we can split any $V\in{\mathcal V}$ as  $V(x)=\widetilde
V(x)+\widetilde W(x)$ where 
\begin{equation}\label{eq:37}
\widetilde
V(x)=\sum_{i=1}^k\frac{\lambda_i\,\alchi_{B(a_i,\delta)}(x)}{|x-a_i|^2}-
\frac{\lambda_{\infty}\alchi_{\R^N\setminus
    B(0,R_0)}(x)}{|x|^2},\quad\d>0,\quad R_0>0,\quad 
\mu(\widetilde V)>0,
\end{equation}
and $\widetilde W\in
L^{N/2}(\R^N)\cap L^{\infty}(\R^N)$.

\begin{Lemma}{\bf [\,Self-adjointness  criterion in ${\mathcal
      V}$\,]}\label{l:sac} 
Let $V\in{\mathcal V}$ and $V=\widetilde V+\widetilde W$, with
$\widetilde V$ as in~(\ref{eq:37}) and $\widetilde W\in
L^{N/2}(\R^N)\cap L^{\infty}(\R^N)$. Then the operator $-\D-V$ is
essentially  
self-adjoint in $C^{\infty}_{\rm c}(\R^N\setminus\{a_1,\dots,a_k\})$
if and only if $\mathop{\rm Range}(-\D-\widetilde V+b)$ is dense in
$L^2(\R^N)$ for some $b>0$.
\end{Lemma}
\begin{pf}
For any $b>0$,  we can
split the operator $-\D-V$ 
as $(-\D-\widetilde V+b)-(\widetilde W+b)$, i.e. as a bounded
perturbation of the positive operator $-\D-\widetilde V+b$. In view of
the Kato-Rellich Theorem (see e.g. \cite[Theorem 4.4]{kato}),
the operator $-\D-V$ is essentially self-adjoint
in $C^{\infty}_{\rm c}(\R^N\setminus\{a_1,\dots,a_k\})$ if and only if
$-\D-\widetilde V+b$ is essentially self-adjoint for some $b>0$.
The conclusion now follows from well-known self-adjointness criteria
for positive operators (see \cite[Theorem X.26]{reedsimon}). 
\end{pf}

\noindent The above criterion provides the following {\em{non 
self-adjointness  condition}} in ${\mathcal V}$.  
\begin{Corollary}\label{l:nsac} 
Let $V\in{\mathcal V}$ and $V=\widetilde V+\widetilde W$, with
$\widetilde V$ as in (\ref{eq:37}) and $\widetilde W\in
L^{N/2}(\R^N)\cap L^{\infty}(\R^N)$. Assume that there exist $v\in
L^2(\R^N)$, $v(x)\geq0$ a.e. in $\R^N$, $\int_{\R^N}v^2>0$, a distribution 
$h\in H^{-1}(\R^N)$, and $b>0$ such that 
\begin{equation}\label{eq:47}
{}_{H^{-1}(\R^N)}\big\langle
h,u\big\rangle_{H^{1}(\R^N)}\leq 0\quad\text{for all }u\in H^1(\R^N),
\ u\geq 0\ \text{a.e in }\R^N,
\end{equation}
and
\begin{equation}\label{eq:46}
-\D v-\tilde V\,v+b \, v=h\quad\text{in }{\mathcal D}'(\R^N\setminus\{a_1,\dots,a_k\}).
\end{equation}
Then the operator $-\D-V$ is not
essentially  
self-adjoint in $C^{\infty}_{\rm c}(\R^N\setminus\{a_1,\dots,a_k\})$.
\end{Corollary}
\begin{pf}
From Lemma \ref{l:sac} and the Kato-Rellich Theorem, it is enough to prove that
$\mathop{\rm Range}(-\D-\widetilde V+b)$ is not dense in $L^2(\R^N)$. To this aim
we will show that $v$ does not belong to the closure of $\mathop{\rm
Range}(-\D-\widetilde V+b)$ in $L^2(\R^N)$. Arguing by
contradiction, we assume that there exist  sequences
$\{u_n\}_n\subset C^{\infty}_{\rm c}(\R^N\setminus\{a_1,\dots,a_k\})$
and $\{f_n\}_n\subset L^2(\R^N)$ such that $f_n\to v$ in $L^2(\R^N)$ and
\begin{equation}\label{eq:45}
-\D u_n(x)-\tilde V(x)\, u_n(x)+b\, u_n(x)=f_n(x).
\end{equation}
In view of the Lax-Milgram Theorem there exists $u\in H^1(\R^N)$, weakly solving 
\begin{equation}\label{eq:44}
-\Delta u(x)-\tilde V(x)u(x)+b\, u(x)=v(x).
\end{equation}
Testing (\ref{eq:44}) with $-u^-$, we easily obtain that $u\geq 0$
a.e. in $\R^N$, hence by the Strong Maximum Principle we deduce that
$u>0$ in $\R^N\setminus \{a_1,\dots,a_k\}$. Subtracting (\ref{eq:44})
from (\ref{eq:45}) and multiplying by $u_n-u$, we find that
$\|u_n-u\|_{H^1(\R^N)}\leq{\rm const\,}\|f_n-v\|_{L^2(\R^N)}$, hence
$u_n\to u$ in $H^1(\R^N)$. Testing~(\ref{eq:46}) with $u_n$ and using (\ref{eq:45}),
we obtain 
$$ 
{}_{H^{-1}(\R^N)}\big\langle
h,u_n\big\rangle_{H^{1}(\R^N)}=
\int_{\R^N}f_n(x)v(x)\,dx,
$$
which, passing to the limit, yields
$$ 
{}_{H^{-1}(\R^N)}\big\langle
h,u\big\rangle_{H^{1}(\R^N)}=
\int_{\R^N}v^2(x)\,dx>0.
$$
The above identity contradicts assumption (\ref{eq:47}).
\end{pf}

\noindent We now extend Theorem \ref{t:kwss} to our class 
of multi-polar potentials, thus proving Theorem \ref{t:self-adjo}.

\bigskip\noindent

\begin{pfn}{Theorem \ref{t:self-adjo}}

\medskip\noindent
{\bf Step 1: } if  $\lambda_i< (N-2)^2/4-1$ for all $i=1,\dots,k$, then
$-\D-V$  is essentially self-adjoint. 

\smallskip\noindent 
In view of
Lemma \ref{l:sac}, to prove essential self-adjointness it is enough to  
show that $\mathop{\rm Range}(-\D-\widetilde V+b)$ is dense in
$L^2(\R^N)$ for some $b>0$, where $\widetilde V$ is as in
(\ref{eq:37}). Let $f\in C^{\infty}_{\rm
  c}(\R^N\setminus\{a_1,\dots,a_k\})$ and $b>0$. By the Lax-Milgram
Theorem, there exists $u\in H^1(\R^N)$ weakly solving  
$$
-\D u(x)-\widetilde V(x)u(x)+b \,u(x)=f(x) \quad\text{in }\R^N.
$$
From Lemma \ref{l:1} we deduce the following asymptotic behavior of
$u$ at poles  
\begin{align}\label{eq:38}
 u(x)\sim
|x-a_i|^{-a_{\lambda_i}},\quad \text{ as }x\to a_i.
\end{align}
Hence the function $g(x):=\widetilde V(x)u(x)-b \,u(x)+f(x)\sim
|x-a_i|^{-a_{\lambda_i}-2}$ as $x\to a_i$. In particular, if $\lambda_i<
(N-2)^2/4-1$ for all $i=1,\dots,k$, then $g\in
L^2(\R^N)$. 
Green's representation formula yields
\begin{align}\label{eq:43}
u(x)=&\frac{1}{N(N-2)\omega_N}\bigg[\int_{B(a_i,\delta)}
\frac{g(y)}{|x-y|^{N-2}}\,dy
+\int_{\partial
  B(a_i,\delta)}\frac{1}{|x-y|^{N-2}}\,\frac{\partial u}{\partial
  \nu}\,ds\bigg]\\
\notag&+\frac{1}{N\omega_N}\int_{\partial
  B(a_i,\delta)}\frac{u(y)}{|x-y|^{N-1}}\,ds,
\quad x\in B(a_i,\delta),
\end{align}
where $\omega_N$ denotes the volume of the unit ball in $\R^N$, 
$\nu$ is the unit outward normal to $\partial
  B(a_i,\delta)$, and $ds$ indicates the $(N-1)$-dimensional area
  element in $\partial   B(a_i,\delta)$.
It is easy to verify that the functions
\begin{align*}
x\mapsto \int_{\partial
  B(a_i,\delta)}\frac{1}{|x-y|^{N-2}}\,\frac{\partial u}{\partial
  \nu}\,ds(y),\qquad
 x\mapsto \int_{\partial
  B(a_i,\delta)}\frac{u(y)}{|x-y|^{N-1}}\,ds(y),
\end{align*}
are of class $C^1( B(a_i,\delta))$. From Lemma \ref{l:green} of the
appendix, we have that
$$
\n\bigg(
\frac{1}{N(N-2)\omega_N}\int_{B(a_i,\delta)} 
\frac{g(y)}{|x-y|^{N-2}}\,dy
\bigg)=-\frac{1}{N\omega_N}\int_{B(a_i,\delta)}
\frac{x-y}{|x-y|^N}g(y)\,dy.
$$
Consequently  
\begin{align}\label{eq:40}
\bigg|\n\bigg(
\frac{1}{N(N-2)\omega_N}\int_{B(a_i,\delta)} 
\frac{g(y)}{|x-y|^{N-2}}\,dy 
\bigg)\bigg|\leq {\rm const\,}\int_{B(a_i,\delta)}\frac{|y-a_i|^{-a_{\lambda_i}-2}}{|x-y|^{N-1}}
\,dy.
\end{align}
If $\lambda_i>1-N$, i.e. $a_{\lambda_i}>-1$, then 
$$
\bigg|\n\bigg(
\frac{1}{N(N-2)\omega_N}\int_{B(a_i,\delta)} 
\frac{g(y)}{|x-y|^{N-2}}\,dy 
\bigg)\bigg|\leq {\rm const\,}h_i(x-a_i),
$$
where $$
h_i(x)=\int_{\R^N}\frac{|y|^{-a_{\lambda_i}-2}}{|x-y|^{N-1}}
\,dy.
$$
An easy scaling argument shows that $h_i(\a x)=\a^{-a_{\lambda_i}-1}h_i(x)$
for all $\a>0$, hence $h_i(x)=|x|^{-a_{\lambda_i}-1}h_i(e_1)$, where
$e_1=(1,0,\dots,0)\in\R^N$. Then, if $\lambda_i>1-N$,
\begin{align}\label{eq:41}
 \bigg|\n\bigg(
\frac{1}{N(N-2)\omega_N}\int_{B(a_i,\delta)} 
\frac{g(y)}{|x-y|^{N-2}}\,dy 
\bigg)\bigg|\leq {\rm const\,} |x-a_i|^{-a_{\lambda_i}-1}.
\end{align}
If $\lambda_i\leq 1-N$, i.e. $a_{\lambda_i}\leq -1$, we fix $0<\e<\frac{N-2}2$
and notice that, from (\ref{eq:40}),  
$$
\bigg|\n\bigg(
\frac{1}{N(N-2)\omega_N}\int_{B(a_i,\delta)} 
\frac{g(y)}{|x-y|^{N-2}}\,dy 
\bigg)\bigg|\leq \delta^{-a_{\lambda_i}-1+\e}\,k_i(x),
$$
where 
$$
k_i(x)=
\int_{\R^N} \frac1{|y-a_i|^{1+\e}|y-x|^{N-1}}\,dy.
$$
An easy scaling argument shows that $k_i(\a x)=\a^{-\e}k_i(x)$
for all $\a>0$, hence $k_i(x)=|x|^{-\e}k_i(e_1)$. Then, if $\lambda_i\leq 1-N$
\begin{align}\label{eq:42}
 \bigg|\n\bigg(
\frac{1}{N(N-2)\omega_N}\int_{B(a_i,\delta)} 
\frac{g(y)}{|x-y|^{N-2}}\,dy 
\bigg)\bigg|\leq C(\e) |x-a_i|^{-\e},
\end{align}
for some positive constant $C(\e)$ depending on $\e$ (and also on $N$,
$\lambda_i$, and
$u$). Representation (\ref{eq:43}), regularity of the boundary terms, and estimates (\ref{eq:41}--\ref{eq:42}) yield
\begin{equation}\label{eq:39}
\n u(x)
=
\begin{cases}
O\big(|x-a_i|^{-a_{\lambda_i}-1}\big),&\text{if }\lambda_i>1-N,\\[5pt]
O\big(|x-a_i|^{-\e}\big),&\text{if }\lambda_i\leq 1-N,
\end{cases}
\qquad \text{ as }x\to a_i.
\end{equation}
For all $n\in\N$ let $\eta_n$ be a cut-off function such
that $\eta_n\in  C^{\infty}_{\rm c}(\R^N\setminus\{a_1,\dots,a_k\})$,
$0\leq\eta_n\leq 1$, and
\begin{align*}
&\eta_n(x)\equiv 0\text{ in }
\bigcup_{i=1}^kB\Big(a_i, \frac1{2n}\Big)\cup \big(\R^N\setminus B(0,2n)\big),\quad
\eta_n(x)\equiv 1\text{ in }B(0,n)\setminus \bigcup_{i=1}^k B\Big(a_i,
\frac1{n}\Big),\\
&|\n \eta_n(x)|\leq C\,n\text{ in } \bigcup_{i=1}^k \Big(B\Big(a_i,
\frac1{n}\Big)\setminus B\Big(a_i, \frac1{2n}\Big)\Big),\quad |\n
\eta_n(x)|\leq \frac{C}n \text{ in }B(0,2n)\setminus B(0,n),\\
&|\D \eta_n(x)|\leq C\,n^2\text{ in } \bigcup_{i=1}^k \Big(B\Big(a_i,
\frac1{n}\Big)\setminus B\Big(a_i, \frac1{2n}\Big)\Big),\quad |\D
\eta_n(x)|\leq \frac{ C}{n^2} \text{ in }B(0,2n)\setminus B(0,n),
\end{align*}
for some positive constant $C$ independent of $n$.
Let us set $f_n:=\eta_n f-2\n \eta_n\cdot\n u-u\,\D\eta_n$ and $u_n:=\eta_n u$,
so that $u_n\in C^{\infty}_{\rm c}(\R^N\setminus\{a_1,\dots,a_k\})$
and  
$
-\D u_n(x)-\widetilde V(x)u_n(x)+b\, u_n(x)= f_n(x)
$.
In particular $f_n\in \mathop{\rm Range}(-\D-\widetilde
V+b)$. Furthermore $\eta_n f\to f$ in $L^2(\R^N)$, while (\ref{eq:38})
and (\ref{eq:39}) yield
\begin{align*}
\int_{\R^N}&|\n \eta_n(x)|^2|\n u(x)|^2\,dx\\
&\leq {\rm const\,} n^2\sum_{i=1}^k\int_{B(a_i,
\frac1{n})\setminus B(a_i, \frac1{2n})}|\n u(x)|^2\,dx+\frac{\rm const\,}{n^2}
\int_{B(0,2n)\setminus B(0,n)}|\n u(x)|^2\,dx\\
&\leq  {\rm const\,}
n^2\bigg[\sum_{\lambda_i>1-N}\int_{B(0,\frac1n)}|x|^{-2a_{\lambda_i}-2}\,dx
+\sum_{\lambda_i\leq 1-N}\int_{B(0,\frac1n)}|x|^{-2\e}\,dx\bigg]
+\frac{\rm
  const\,}{n^2}\|u\|_{H^1(\R^N)}\\
&\leq  {\rm const\,}\bigg[\sum_{\lambda_i>1-N } n^{2
  a_{\lambda_i}+4-N}+\sum_{\lambda_i\leq 1-N } n^{2
  +2\e-N}
+n^{-2}\bigg]
\end{align*}
and
\begin{align*}
\int_{\R^N}&|\D \eta_n(x)|^2|u(x)|^2\,dx\\
&\leq {\rm const\,} n^4\sum_{i=1}^k\int_{B(a_i,
\frac1{n})\setminus B(a_i, \frac1{2n})}|u(x)|^2\,dx+\frac{\rm const\,}{n^4}
\int_{B(0,2n)\setminus B(0,n)}|u(x)|^2\,dx\\
&\leq  {\rm const\,}
n^4\sum_{i=1}^k\int_{B(0,1/n)}|x|^{-2a_{\lambda_i}}\,dx+\frac{\rm
  const\,}{n^4}\|u\|_{H^1(\R^N)}\\
&\leq  {\rm const\,}\bigg[\sum_{i=1}^k n^{2
  a_{\lambda_i}+4-N}+n^{-4}\bigg].
\end{align*}
Since for $\lambda_i< (N-2)^2/4-1$ there holds $2 a_{\lambda_i}+4-N<0$, we
conclude that $f_n\to f$ in $L^2(\R^N)$. Hence $\mathop{\rm Range}(-\D-\widetilde
V+b)$ is dense in $C^{\infty}_{\rm
  c}(\R^N\setminus\{a_1,\dots,a_k\})$. Since $C^{\infty}_{\rm
  c}(\R^N\setminus\{a_1,\dots,a_k\})$ is dense in $L^2(\R^N)$, we
obtain that  $\mathop{\rm Range}(-\D-\widetilde
V+b)$ is dense in $L^2(\R^N)$.

\medskip\noindent {\bf Step 2: } if $\lambda_i\leq (N-2)^2/4-1$ for all
$i\in\{1,\dots,k\}$, then $-\D-V$  is essentially self-adjoint.

\smallskip\noindent
Let us fix $b>0$, $f\in C^{\infty}_{\rm
  c}(\R^N\setminus\{a_1,\dots,a_k\})$. To prove essential
self-adjointness it is enough to find some $g\in \mathop{\rm
  Range}(-\D-\widetilde V+b)$ such that $g$ is arbitrarily closed to
$f$ in $L^2(\R^N)$. To
this aim, we fix $\e>0$ and  notice that there exists $0<\sigma<1$
such that if  $u\in
H^1(\R^N)$ solves 
\begin{equation}\label{eq:54}
-\D u(x)-\widetilde V_{\sigma}(x)u(x)+b\,u(x)=f(x), 
\end{equation} 
where $$
\widetilde V_{\sigma}(x):=
\sum_{\lambda_i=\left(\!\frac{N-2}2\!\right)^2-1}\frac{(\lambda_i-\sigma)\, 
\alchi_{B(a_i,\delta)}(x)}{|x-a_i|^2}+\sum_{\lambda_i<
\left(\!\frac{N-2}2\!\right)^2-1}
\frac{\lambda_i  \,\alchi_{B(a_i,\delta)}(x)}{|x-a_i|^2}-
\frac{\lambda_{\infty}\alchi_{\R^N\setminus
    B(0,R_0/\delta)}(x)}{|x|^2},
$$
then 
\begin{equation}\label{eq:53}
\|(\widetilde V_{\sigma}-\widetilde V)u\|_{L^2(\R^N)}<\e.
\end{equation}
Indeed, by Remark \ref{r:h1}, there exists a positive constant $C$
independent on $\sigma\in(0,1)$, such that all solutions of (\ref{eq:54})
can be estimated as 
$$
|u(x)|\leq C\,|x-a_i|^{-a_{(\lambda_i-\sigma)}}\|u\|_{H^1(\R^N)}\quad
\text{in }B(a_i,\delta),
$$
for all $i$ such that
$\lambda_i=\left(\!\frac{N-2}2\!\right)^2-1$. Moreover, testing
(\ref{eq:54}) by $u$ there results that all solutions of
(\ref{eq:54}) satisfy
$$
\|u\|_{H^1(\R^N)}\leq\frac{\|f\|_{L^2(\R^N)}}{\min\{\mu(\widetilde
  V),b\}}.
$$
 Then
for all $i$ such that $\lambda_i=\left(\!\frac{N-2}2\!\right)^2-1$, we have
that
\begin{align*}
\left\| \frac{\sigma
  \,  \alchi_{B(a_i,\delta)}(x)\,u}{|x-a_i|^2}\right\|_{L^2(\R^N)}^2
&\leq C^2\sigma^2\|u\|_{H^1(\R^N)}^2\int_0^1
r^{N-5-2a_{(\lambda_i-\sigma)}}\,dr\\
&=\frac{C^2\,\|f\|^2_{L^2(\R^N)}}{2\big(\min\{\mu(\widetilde
  V),b\}\big)^2}\,\frac{\sigma^2}{\sqrt{1+\sigma}-1}=o(1)\quad\text{as }\sigma\to0.
\end{align*}
Therefore it is possible to choose $\sigma$ small enough in order to
ensure that all solutions of (\ref{eq:54}) satisfy~(\ref{eq:53}).
For such a $\sigma$, let $u\in H^1(\R^N)$ be a solution to
(\ref{eq:54}).
Let $\eta_n$ be the sequence of cut-off functions introduced in step 1. As
in step 1, we have that  
$f_n:=\eta_n f-2\n \eta_n\cdot\n u-u\,\D\eta_n$ converges to $f$ in
$L^2(\R^N)$. Hence, for every $n$ large enough,
$\|f_n-f\|_{L^2(\R^N)}<\e$. 
Moreover $u_n:=\eta_n u\in C^{\infty}_{\rm
  c}(\R^N\setminus\{a_1,\dots,a_k\})$ and
$$\|(\widetilde V_{\sigma}-\widetilde V)u_n\|_{L^2(\R^N)}\leq
\|(\widetilde V_{\sigma}-\widetilde V)u\|_{L^2(\R^N)}<\e$$ 
for any $n$. Setting $g_n(x):=f_n(x)+(\widetilde
V_{\sigma}(x)-\widetilde V(x))u_n(x)$, we obtain that $u_n$ satisfies 
$$
-\D u_n(x)-\widetilde V(x) u_n(x)+b\,u_n(x)=g_n(x),
$$
i.e. $g_n\in  \mathop{\rm
  Range}(-\D-\widetilde V+b)$, and $\|g_n-
f\|_{L^2(\R^N)}<2\e$ for 
large $n$. The proof of step 2 is thereby complete. 

\medskip\noindent {\bf Step 3: } if $\lambda_i> (N-2)^2/4-1$ for some
$i\in\{1,\dots,k\}$, then $-\D-V$  is not essentially self-adjoint.

\smallskip\noindent
Let $V=\widetilde V+\widetilde W$, with
$\widetilde V$ as in (\ref{eq:37}) and $\widetilde W\in
L^{N/2}(\R^N)\cap L^{\infty}(\R^N)$. Let us fix $i\in\{1,\dots,k\}$ such
that  $\lambda_i> (N-2)^2/4-1$ and $\a<0$, and consider the solution  $\psi\in
C^1\big((-\infty,\ln\d]\big)$ of the Cauchy problem 
$$
\begin{cases}
\psi''(s)-\omega_{\lambda_i}^2\,\psi(s)=b\, e^{2s}\,\psi(s),\\
\psi(\ln \d)=0,\quad
\psi'(\ln \d)=\a,
\end{cases}
$$
where $\omega_{\lambda_i}:=\sqrt{\big(\frac{N-2}2\big)^2-\lambda_i}$ and $\d$ is
as in (\ref{eq:37}). In view of Lemma \ref{l:ode} of the appendix, we can 
estimate $\psi$ as 
\begin{equation}\label{eq:48}
0\leq \psi(s)\leq C\,
\,e^{-\omega_{\lambda_i} s}\quad\text{for all }s\leq\ln\d,
\end{equation}
for some positive constant $C=C(\lambda_i,\d,\a,b)$. Let us set
$$
v(x):=
\begin{cases}
|x-a_i|^{-\frac{N-2}2}\,\psi(\ln|x-a_i|),&\text{if }x\in B(a_i,\d)\setminus\{a_i\},\\
0,&\text{if }x\in \R^N\setminus\overline{B(a_i,\d)}.
\end{cases}
$$
From (\ref{eq:48}) we infer that 
\begin{equation}\label{eq:49}
0\leq v(x)\leq
C|x-a_i|^{-\frac{N-2}2-\sqrt{\left(\frac{N-2}2\right)^2-\lambda_i}}\quad\text{in }
B(a_i,\d). 
\end{equation}
The assumption  $\lambda_i> (N-2)^2/4-1$ and estimate (\ref{eq:49}) ensure
that $v\in L^2(\R^N)$. Moreover the restriction of $v$ to $B(a_i,\d)$
satisfies 
$$
\begin{cases}
-\D v(x)-
\dfrac{\lambda_i}{|x-a_i|^2}\,v(x)+b\,v(x)=0,&\text{in } B(a_i,\d),\\[10pt]
v=0\quad\text{and}\quad \dfrac{\partial v}{\partial\nu}=\d^{-\frac
  N2}\a,&\text{on } \partial B(a_i,\d).
\end{cases}
$$
As a consequence the distribution $-\D v-\widetilde
V\,v+b\,v\in{\mathcal D}'(\R^N\setminus\{a_1,\dots,a_k\})$ acts as follows: 
$$ 
{}_{{\mathcal D}'(\R^N\setminus\{a_1,\dots,a_k\})}\big\langle
-\D v-\widetilde
V\,v+b\,v\,,\,\varphi\big\rangle_{C^{\infty}_{\rm
  c}(\R^N\setminus\{a_1,\dots,a_k\})}=\d^{-\frac
  N2}\a \int_{\partial B(a_i,\d)}\varphi(x)\,ds.
$$
Hence $h=-\D v-\widetilde
V\,v+b\,v\in H^{-1}(\R^N)$ and  satisfies
(\ref{eq:47}) as $\a<0$. From Corollary \ref{l:nsac}, we finally deduce that 
the operator $-\D-V$ is not
essentially  
self-adjoint in $C^{\infty}_{\rm c}(\R^N\setminus\{a_1,\dots,a_k\})$.
\end{pfn}

The following theorem characterizes essential self-adjointness of 
Schr\"odinger operators with potentials carrying infinitely many
singularities distributed on reticular structures.  

\begin{Theorem}\label{t:self-adjo-ret}
For $\lambda<(N-2)^2/4$ and $\{a_n\}_n\subset\R^N$ satisfying 
(\ref{eq:55}) and $|a_n-a_m|\geq 1$
for all $n\neq m$, let $0<\delta<1/2$ be given by Lemma \ref{l:ret} and  
$$
V(x)=\lambda\sum_{n=1}^{\infty}\frac{\alchi_{B(a_n,\delta)}(x)}{|x-a_n|^2}.
$$
Then $-\D-V$  is essentially self-adjoint
in $C^{\infty}_{\rm c}\left(\R^N\setminus\{a_n\}_{n\in\N}\right)$ if
and only if $\lambda\leq (N-2)^2/4-1$.
\end{Theorem}
\begin{pf}
From the Kato-Rellich Theorem 
the operator $-\D-V$ is essentially self-adjoint
in $C^{\infty}_{\rm c}\left(\R^N\setminus\{a_n\}_{n\in\N}\right)$ if
and only if 
$-\D- V+b$ is essentially self-adjoint for some $b>0$.
In view of Lemma \ref{l:ret}, for any $b>0$, $-\D- V+b$ is
positive. Hence essential self-adjointness is equivalent to density of
$\mathop{\rm Range}(-\D- V+b)$ in $L^2(\R^N)$ for some $b>0$.

Let us first prove that,  if  $\lambda< (N-2)^2/4-1$, then
$-\D-V$  is essentially self-adjoint.  For $f\in C^{\infty}_{\rm
  c}\left(\R^N\setminus\{a_n\}_{n\in\N}\right)$ and $b>0$, the
Lax-Milgram Theorem provides a unique $u\in
H^1(\R^N)$  weakly solving
$$
-\D u(x)-V(x)u(x)+b \,u(x)=f(x) \quad\text{in }\R^N.
$$
From Lemma \ref{l:1} and arguing as in the proof of Theorem
\ref{t:self-adjo}, we deduce that
\begin{align}\label{eq:100}
 u(x)\sim
|x-a_n|^{-a_{\lambda}},\quad\text{and}\quad
\n u(x)
=
\begin{cases}
O\big(|x-a_n|^{-a_{\lambda}-1}\big),&\text{if }\lambda>1-N,\\[5pt]
O\big(|x-a_n|^{-\e}\big),&\text{if }\lambda\leq 1-N,
\end{cases}
\qquad \text{ as }x\to a_n,
\end{align}
where $0<\e<\frac{N-2}{2}$. Since $\lambda<(N-2)^2/4-1$, we have that $2a_{\lambda}+4-N<0$, 
hence, for all $j\in\N$, we can choose $N_j\in\N$ such that
$N_j\to+\infty$ as 
$j\to\infty$, $N_jj^{2a_{\lambda}+4-N}\to 0$, and $N_j j^{2\e-N+2}\to
0$, and let $R_j>0$ such that $R_j\to+\infty$ as $j\to\infty$ and $B(a_n,1/j)\subset B(0,R_j)$ for all
$n=1,\dots, N_j$.  Let $\eta_j$ be a cut-off function such 
that $\eta_j\in  C^{\infty}_{\rm c}\left(\R^N\setminus\{a_n\}_{n\in\N}\right)$,
$0\leq\eta_j\leq 1$, and
\begin{align*}
&\eta_j(x)\equiv 0\text{ in }
\bigcup_{n=1}^{N_j}B\Big(a_n, \frac1{2j}\Big)\cup (\R^N\setminus
B(0,2R_j)),\quad 
\eta_j(x)\equiv 1\text{ in }B(0,R_j)\setminus \bigcup_{n=1}^{N_j} B\Big(a_n,
\frac1{j}\Big),\\
&|\n \eta_j(x)|\leq C\,j\text{ in } \bigcup_{n=1}^{N_j} \Big(B\Big(a_n,
\frac1{j}\Big)\setminus B\Big(a_n, \frac1{2j}\Big)\Big),\quad
|\n \eta_j(x)|\leq \frac{C}{R_j}\text{ in }B(0,2 R_j) \setminus B(0,R_j),
\\
&|\D \eta_j(x)|\leq C\,j^2\text{ in } \bigcup_{n=1}^{N_j} \Big(B\Big(a_n,
\frac1{j}\Big)\setminus B\Big(a_n, \frac1{2j}\Big)\Big),\quad
|\D \eta_j(x)|\leq \frac{C}{R_j^2}\text{ in }B(0,2 R_j) \setminus B(0,R_j),
\end{align*}
for some positive constant $C$ independent of $j$ and $n$.
Let $f_j:=\eta_j f-2\n \eta_j\cdot\n u-u\,\D\eta_j$ and $u_j:=\eta_j u$,
so that $u_j\in C^{\infty}_{\rm c}\left(\R^N\setminus\{a_n\}_{n\in\N}\right)$
and  
$
-\D u_j(x)- V(x)u_j(x)+b\, u_j(x)= f_j(x)
$.
In particular $f_j\in \mathop{\rm Range}(-\D-V+b)$. Furthermore $\eta_j f\to f$ in $L^2(\R^N)$, while (\ref{eq:100}) 
yields
\begin{align*}
\int_{\R^N}&|\n \eta_j(x)|^2|\n u(x)|^2\,dx\\
&\leq {\rm const\,} j^2\sum_{n=1}^{N_j}\int_{B(a_n,
\frac1{j})\setminus B(a_n, \frac1{2j})}|\n u(x)|^2\,dx+\frac{\rm
const\,}{R_j^2} 
\int_{B(0,2R_j)\setminus B(0,R_j)}|\n u(x)|^2\,dx
\\[5pt]
&\leq 
\begin{cases}
 {\rm const\,}\big[
N_j\, j^{2a_{\lambda}+4-N}+R_j^{-2}\|u\|_{H^1(\R^N)}\big],&\text{ if }\lambda>1-N\\[3pt]
{\rm const\,}\big[N_j\, j^{2\e-N+2}+R_j^{-2}\|u\|_{H^1(\R^N)}\big],&\text{ if }\lambda\leq 1-N
\end{cases}
\end{align*}
and, in a similar way,
\begin{align*}
\int_{\R^N}|\D \eta_j(x)|^2|u(x)|^2\,dx\leq  {\rm const\,}\big[N_j\, j^{2a_{\lambda}+4-N}+R_j^{-4}\big].
\end{align*}
By the choice of $N_j$, we deduce that $f_j\to f$ in $L^2(\R^N)$. 
Hence $\mathop{\rm Range}(-\D-\
V+b)$ is dense in $C^{\infty}_{\rm
  c}\left(\R^N\setminus\{a_n\}_{n\in\N}\right)$ and consequently  in $L^2(\R^N)$.

\medskip\noindent To prove essential self-adjointness for $\lambda=
(N-2)^2/4-1$, we can argue as in the proof of Theorem
\ref{t:self-adjo}, Step 2, i.e. by approximation of the resonant potential $V$
 with sub-resonant potentials. To do that, we need to prove that  
for fixed $b>0$, $f\in C^{\infty}_{\rm
  c}\left(\R^N\setminus\{a_n\}_{n\in\N}\right)$, and $\e>0$, there
exists $0<\sigma<1$  such that if $u\in H^1(\R^N)$ solves
\begin{equation}\label{eq:57}
-\D u(x)- V_{\sigma}(x)u(x)+b\,u(x)=f(x),\quad\text{where }
V_{\sigma}(x):=(\lambda-\sigma)\sum_{n=1}^{\infty}\frac{  
\alchi_{B(a_n,\delta)}(x)}{|x-a_n|^2},
\end{equation} 
then 
\begin{equation}\label{eq:58}
\|(V_{\sigma}- V)u\|_{L^2(\R^N)}<\e.
\end{equation}
Indeed, by Remark \ref{r:h1}, there exists a positive constant $C$
independent on $\sigma\in(0,1)$ and $n\in\N$, such that all solutions
of (\ref{eq:57}) can be estimated as 
$$
|u(x)|\leq C\,|x-a_n|^{-a_{(\lambda-\sigma)}}\|u\|_{H^1(B(a_n,\d'))}\quad
\text{in }B(a_n,\delta),
$$
for some $\delta<\delta'<1/2$ and for all $n\in\N$. Consequently
\begin{align*}
\left\| \frac{\sigma
  \,  \alchi_{B(a_n,\delta)}(x)\,u}{|x-a_n|^2}\right\|_{L^2(\R^N)}^2
&\leq
\frac{C^2}2\|u\|_{H^1(B(a_n,\d'))}^2\,\frac{\sigma^2}{\sqrt{1+\sigma}-1},
\end{align*}
and hence
\begin{align*}
\|(V_{\sigma}- V)u\|_{L^2(\R^N)}&\leq
\frac{C}{\sqrt2}\frac{\sigma}{\sqrt{\sqrt{1+\sigma}-1}}\|u\|_{H^1(\R^N)}
\\
&\leq {\rm
  const\,}\|f\|_{L^2(\R^N)}\frac{\sigma}{\sqrt{\sqrt{1+\sigma}-1}}
\to 0\quad\text{as }\sigma\to0. 
\end{align*}
Therefore it is possible to choose $\sigma$ small enough in order to
ensure that all solutions of (\ref{eq:57}) satisfy~(\ref{eq:58}).
In order to prove self-adjointness, it is now sufficient to repeat the
argument of Theorem
\ref{t:self-adjo}, Step 2.

\medskip\noindent The proof of non essential self-adjointness in the
case $\lambda> (N-2)^2/4-1$ can be obtained just by mimicking the arguments
of the proof of Theorem \ref{t:self-adjo} and Corollary \ref{l:nsac}.\end{pf}

\section{Spectrum of Schr\"odinger operators with potentials in
${\mathcal V}$}\label{sec:spectr-schr-oper-1}

In this section we study the spectrum of the {\em {Friedrichs
extension}} $(-\D-V)^F$ of Schr\"odinger operators with potentials in
  ${\mathcal V}$, see  \eqref{eq:52}. We recall that in view of
  Theorem \ref{t:self-adjo}, if  $\lambda_i\leq(N-2)^2/4-1$ for all
$i=1,\dots,k$,  then $(-\D-V)^F$ is the only self-adjoint
extension of $-\D-V$. On the other hand,  if $\lambda_i>(N-2)^2/4-1$ for some $i$, then $-\D-V$ has  many self-adjoint extensions,
among which the Friedrichs extension is the only one with domain included
in $H^1(\R^N)$. Due to self-adjointness, the spectrum of 
$(-\D-V)^F$ turns out to be a subset of $\R$, which will
be described below.

\subsection{Essential spectrum}\label{sec:spectr-schr-oper}

Let us start by studying the essential spectrum of the Friedrichs
extension of operators $-\D-V$,
$V\in {\mathcal V}$, in $L^2(\R^N)$. 
The Friedrichs extension $(-\D-V)^F:\ D\big((-\D-V)^F\big)\to
L^2(\R^N)$ defined in~(\ref{eq:52}) is self-adjoint. 
As a consequence, the essential spectrum $\sigma_{\rm
  ess}\big((-\D-V)^F\big)$ can be 
characterized by terms of the {\em Weyl sequences} as follows:
$\lambda\in \sigma_{\rm ess}\big((-\D-V)^F\big)$ if and only if 
\begin{equation}\label{eq:24}
\left\{
\hskip-10pt
\begin{array}{ll}
&\text{there exists
}\{f_n\}_n\subset D\big((-\D-V)^F\big)\text{ such that }\liminf_{n\to+\infty}
\|f_n\|_{L^2(\R^N)}>0,\\[5pt]
& f_n\weakly 0 \text{ weakly in }L^2(\R^N), 
\text{ and }\|-\D f_n-V f_n-\lambda f_n\|_{L^2(\R^N)}\to0,
\end{array}
\right.
\end{equation}
see \cite[p. 167]{davies}.

\medskip\noindent
\begin{pfn}{Proposition \ref{p:spe}.1}

\medskip\noindent{\bf Step 1:} $[0,+\infty)\subseteq\sigma_{\rm ess}(-\D-V)$.
Let $\lambda\geq 0$. It is well known that $\sigma_{\rm ess}(-\D)=[0,+\infty)$, where
 $$-\D:\ D(-\D)=H^2(\R^N)\to L^2(\R^N).
$$
  Hence $\lambda\in \sigma_{\rm ess}(-\D)$ and 
 the characterization given in (\ref{eq:24}) yields a sequence 
$\{f_n\}_n\subset H^2(\R^N)$, such that $\|f_n\|_{L^2(\R^N)}=1$,
$f_n\weakly 0$  weakly in $L^2(\R^N)$ and $\|-\D f_n-\lambda
f_n\|_{L^2(\R^N)}\to0$. By density of $C^{\infty}_{\rm c}(\R^N)$ in
$H^2(\R^N)$, for any $n$ there exists $g_n\in  C^{\infty}_{\rm
  c}(\R^N)$ such that $\|g_n-f_n\|_{H^2(\R^N)}\leq 1/n$.
It is easy to verify that $g_n\weakly 0$  weakly in $L^2(\R^N)$, 
$1/2\leq\|g_n\|_{L^2(\R^N)}\leq 2$ for sufficiently large $n$, and
$\|-\D g_n-\lambda
g_n\|_{L^2(\R^N)}\to0$.
Let us choose a sequence $\{x_n\}_n\subset \R^N$ such that 
\begin{equation}\label{eq:28}
\mathop{\rm supp}\varphi_n\subset \R^N\setminus
B(0,n),\quad\text{where}\quad \varphi_n(x):=g_n(x+x_n).
\end{equation}
By (\ref{eq:28}), it is  easy to prove that $\varphi_n\weakly 0$
weakly in $L^2(\R^N)$ and 
 $1/2\leq\|\varphi_n\|_{L^2(\R^N)}=\|g_n\|_{L^2(\R^N)}\leq 2$.
From (\ref{eq:28}), it follows also that, if $n$ is sufficiently
large, the support of $\varphi_n$ is disjoint 
from all balls $B(a_i,r_i)$ where singularities of $V$ are
located. Therefore 
$$
V\varphi_n=\frac{\lambda_{\infty}\alchi_{\R^N\setminus
  B(0,R)}}{|x|^{2}}\,\varphi_n+W\varphi_n
\in L^2(\R^N)
$$
 and hence
$\varphi_n\in D\big((-\D-V)^F\big)$. 

Furthermore, letting $h_n= -\D \varphi_n-\lambda
\varphi_n$, we have that $\|h_n\|_{L^2(\R^N)}=\|-\D g_n-\lambda
g_n\|_{L^2(\R^N)}\to 0$, hence
$$
\int_{\R^N}|\n
\varphi_n(x)|^2\,dx=
\lambda\int_{\R^N}\varphi_n^2(x)\,dx+\int_{\R^N}h_n(x)\varphi_n(x)\,dx\leq
4\big(\lambda+o(1)\big),\quad\text{as }n\to+\infty.
$$
Since $\varphi_n$ is bounded in $\Di$ and $W\in L^{N/2}(\R^N)\cap
L^{\infty}(\R^N)$, from  Sobolev's inequality we obtain 
\begin{equation}\label{eq:30}
\|W\varphi_n \|_{L^2(\R^N)}^2\leq S^{-1}\|W\|_{L^{\infty}(\R^N)}
\bigg(\int_{\R^N\setminus B(0,n)}|W(x)|^{N/2}\,dx\bigg)^{2/N}\int_{\R^N}|\n
\varphi_n(x)|^2\,dx\to 0
\end{equation}
as $n\to+\infty$. Moreover 
\begin{equation}\label{eq:31}
\bigg\|\frac{\lambda_{\infty}\alchi_{\R^N\setminus
  B(0,R)}}{|x|^{2}}\,\varphi_n
 \bigg\|_{L^2(\R^N)}^2\leq
\frac{\lambda_{\infty}^2}{n^4}\|\varphi_n\|^2_{L^2(\R^N)}\to0\quad\text{as}\quad
n\to+\infty.
\end{equation}
From (\ref{eq:30}--\ref{eq:31}) we deduce that
$\lim_{n\to+\infty}\|V\varphi_n \|_{L^2(\R^N)}=0$. 
As a consequence
\begin{align*}
\|-\D \varphi_n-V\varphi_n-\lambda
\varphi_n\|_{L^2(\R^N)}&\leq \|h_n\|_{L^2(\R^N)}+\|V\varphi_n \|_{L^2(\R^N)}
\to0
\end{align*}
as $n\to+\infty$. Thus, $\{\varphi_n\}_n$ is a Weyl's sequence and
$\lambda\in\sigma_{\rm ess}\big((-\D-V)^F\big)$.

\medskip\noindent{\bf Step 2:} $\sigma_{\rm
  ess}\big((-\D-V)^F\big)\subseteq[0,+\infty)$.
Assume now that $\lambda\in\sigma_{\rm ess}\big((-\D-V)^F\big)$. Then, from
(\ref{eq:24}) there exists a sequence  
$\{f_n\}_n\subset D\big((-\D-V)^F\big)$  such that $\|f_n\|_{L^2(\R^N)}=1$, 
$f_n\weakly 0$ weakly in $L^2(\R^N)$, and $h_n:=-\D f_n-V f_n-\lambda
f_n\to0$ strongly in $L^2(\R^N)$.
By Lemma \ref{l:semibounded}, we can  write $V(x)=\widetilde
V(x)+\widetilde W(x)$ where $\mu(\widetilde V)>0$ and $\widetilde W\in
L^{N/2}(\R^N)\cap L^{\infty}(\R^N)$. Hence
\begin{align*}
\mu(\widetilde V)\int_{\R^N}|\n f_n(x)|^2\,dx&\leq 
\int_{\R^N}\big(|\n f_n(x)|^2-\widetilde V(x)|f_n(x)|^2\big)\,dx\\
&
=\lambda \int_{\R^N}|f_n(x)|^2\,dx+ \int_{\R^N}\widetilde W(x)|f_n(x)|^2\,dx
+\int_{\R^N}h_n(x)f_n(x)\,dx\\
&\leq \lambda+\|\widetilde W\|_{L^{\infty}(\R^N)}+o(1)\quad\text{as }n\to+\infty.
\end{align*}
Being $\{f_n\}_n$ bounded in $H^1(\R^N)$, there exists $f\in
H^1(\R^N)$ such that, up to a subsequence, $f_n\weakly f$ weakly in
$H^1(\R^N)$. Weak convergence of $f_n$ to $0$ in $L^2(\R^N)$ implies
that $f=0$, hence $f_n\weakly 0$ weakly in $H^1(\R^N)$ and a.e. in $\R^N$.
For any measurable set $\omega$, Sobolev's inequality implies
\begin{equation*}
\int_{\omega}\widetilde W(x) f_n^2(x)\,dx
\leq S^{-1}\bigg(\int_{\omega}|\widetilde
W(x)|^{N/2}\,dx\bigg)^{\!\!2/N}\!\!\!\!\int_{\R^N}|\n 
f_n(x)|^2\,dx\leq {\rm const\,}\bigg(\int_{\omega}|\widetilde
W(x)|^{N/2}\,dx\bigg)^{\!\!2/N}\!\!, 
\end{equation*}
hence the integral in left hand side goes to zero both for the Lebesgue
measure of $\omega$ tending to $0$  and for $\omega$ being the
complement of balls with radius tending to $+\infty$. As a consequence,
the Vitali's Convergence Theorem yields 
\begin{equation}\label{eq:33}
\lim_{n\to+\infty}\int_{\R^N}\widetilde W(x) f_n^2(x)\,dx=0.
\end{equation}
From (\ref{eq:33}) and the strong convergence of $h_n$ to $0$ in
$L^2(\R^N)$, we obtain
\begin{align*}
-\lambda\leq \mu(\widetilde V)\int_{\R^N}&|\n f_n(x)|^2\,dx
-\lambda \int_{\R^N}|f_n(x)|^2\,dx\\
&\leq 
\int_{\R^N}\big(|\n f_n(x)|^2-\widetilde V(x)|f_n(x)|^2\big)\,dx
-\lambda \int_{\R^N}|f_n(x)|^2\,dx
\\
&
=\int_{\R^N}\widetilde W(x)|f_n(x)|^2\,dx
+\int_{\R^N}h_n(x)f_n(x)\,dx=o(1)\quad\text{as }n\to+\infty.
\end{align*}
Letting $n\to+\infty$, we obtain $\lambda\geq0$.
\end{pfn}

\begin{remark}{\bf [\,Essential spectrum in the case of infinitely
    many reticular singularities\,]}\label{r:esssperet} 
For $\lambda<(N-2)^2/4$, let $\{a_n\}_n\subset\R^N$ be a sequence of poles
located on a periodic $M$-dimensional 
reticular structure, $M< N-2$. As observed in Remark \ref{rem:ret2}, 
(\ref{eq:55}) is satisfied and Lemma \ref{l:ret} yields
$\d>0$ such that the quadratic form associated to the infinitely
singular operator $-\D-V$,
$V(x)=\lambda\sum_{n=1}^{\infty}|x-a_n|^{-2}\alchi_{B(a_n,\delta)}(x)$,
is positive definite in $\Di$. Since the reticulation does not fill
the whole $\R^N$, we can repeat the translation argument in Step 1 of
the proof of 
Proposition \ref{p:spe}.1  to construct Weyl's sequences. In addition,
the positivity of the quadratic form allows us to mimic the procedure 
developed in Step 2, thus obtaining that the essential spectrum of the
Friedrichs extension $(-\D-V)^F$ is given by the half line $[0,+\infty)$.
\end{remark}

\subsection{Discrete  spectrum}\label{sec:discr-spectr-schr-oper}

\noindent If $\nu_1(V)<0$, then the spectrum of $(-\D-V)^F$ below $0$,
namely the 
{\em discrete spectrum}
$$
\sigma_{\rm d}\big((-\D-V)^F\big):=\sigma\big((-\D-V)^F\big)\setminus
\sigma_{\rm ess}\big((-\D-V)^F\big)=\sigma\big((-\D-V)^F\big)\cap(-\infty,0),
$$
is not empty and
is described as a sequence of eigenvalues
$$
\nu_1(V)<\nu_2(V)<\dots<\nu_k(V)\dots
$$
which admit the following
variational characterization:
\begin{equation*}
\nu_k(V):=\inf_{u\in E^k \setminus\{0\}}\frac{\displaystyle\int_{\R^N}\big(|\n
  u(x)|^2-V(x)u^2(x)\big)\,dx}{\displaystyle\int_{\R^N}|u(x)|^2\,dx},
\quad k=1,2,\dots,
\end{equation*}  
where 
$$
E^k=\bigg\{w\in H^1(\R^N):\ \int_{\R^N}w(x)v_i(x)=0\quad\text{for
}i=1,\dots,k-1\bigg\}
$$
and $\{v_i,\ i=1,\dots,k-1\}$, are the
first $k-1$ eigenfunctions. 
The following corollary of Lemma \ref{l:semibounded} states that whenever
$\nu_1(V)<0$, then it is attained. The corresponding eigenfunction
thus provides  a {\em bound state} in $L^2(\R^N)$ with negative energy.

\begin{Corollary}\label{c:spec}
If $V\in{\mathcal V}$ and
$\nu_1(V)<0$, 
then $\nu_1(V)$ is attained.
\end{Corollary}
\begin{pf} 
In view of Lemma \ref{l:semibounded}, we can write $V$ as  $V(x)=\widetilde
V(x)+\widetilde W(x)$ where $\mu(\widetilde V)>0$ and $\widetilde W\in
L^{N/2}(\R^N)\cap L^{\infty}(\R^N)$. 
Let $\{u_n\}_n\subset H^1(\R^N)$ be a minimizing sequence such that
\begin{align*}
\int_{\R^N}|u_n(x)|^2\,dx=1\quad\text{and}\quad 
\lim_{n\to+\infty}\int_{\R^N}\big(|\n
u_n(x)|^2-V(x)u_n^2(x)\big)\,dx=\nu_1(V).
\end{align*}
Since 
\begin{align*}
\mu(\widetilde V)\int_{\R^N}|\n u_n(x)|^2\,dx&\leq
\int_{\R^N}\big(|\n u_n(x)|^2-\widetilde V(x)u_n^2(x)\big)\,dx\\
&=\nu_1(V)+\int_{\R^N}\widetilde W(x)|u_n(x)|^2\,dx+o(1)\leq
\nu_1(V)+\|\widetilde W\|_{L^{\infty}(\R^N)}+o(1)
\end{align*}
as $n\to+\infty$, we obtain that $\{u_n\}_n$ is bounded in
$H^1(\R^N)$, hence, up to a subsequence, $u_n\weakly u$ weakly in
$H^1(\R^N)$ and a.e. in $\R^N$. Vitali's Convergence Theorem easily yields
$$\int_{\R^N}\widetilde W(x)|u_n(x)|^2\,dx\to \int_{\R^N}\widetilde
W(x)|u(x)|^2\,dx.
$$
Therefore, taking into account that
$\int_{\R^N}\big(|\n u(x)|^2-\widetilde V(x)u^2(x)\big)\,dx$ is an
equivalent norm, we deduce
\begin{multline}\label{eq:2}
\int_{\R^N}\big(|\n  u(x)|^2-\widetilde V(x)u^2(x)\big)\,dx- \int_{\R^N}\widetilde
W(x)|u(x)|^2\,dx\\
\leq \liminf_{n\to+\infty}\bigg(\int_{\R^N}\big(|\n  u_n(x)|^2-\widetilde V(x)u_n^2(x)\big)\,dx\bigg)
-\lim_{n\to+\infty}\int_{\R^N}\widetilde W(x)|u_n(x)|^2\,dx=\nu_1(V)<0.
\end{multline}
Hence $u\not\equiv 0$. Then from (\ref{eq:20}) and (\ref{eq:2}) it
follows
\begin{equation}\label{eq:21}
\nu_1(V)\int_{\R^N}|u(x)|^2\,dx\leq 
\int_{\R^N}\big(|\n  u(x)|^2-V(x)u^2(x)\big)\,dx\leq \nu_1(V).
\end{equation}
Hence $\int_{\R^N}|u(x)|^2\,dx\geq1$. On the other hand, by weakly
lower semi-continuity 
of the $L^2$--norm, we have that $\int_{\R^N}|u(x)|^2\,dx\leq
\liminf_{n\to+\infty}\int_{\R^N}|u_n(x)|^2\,dx=1$. Therefore
$\int_{\R^N}|u(x)|^2\,dx=1$ and, from (\ref{eq:21}),
$$\int_{\R^N}\big(|\n  u(x)|^2-V(x)u^2(x)\big)\,dx=\nu_1(V),$$ 
i.e. $u$ attains the infimum in (\ref{eq:20}). 
\end{pf}

\noindent Fix an integer $k\geq 1$.
Arguing as above, 
Lemma \ref{l:semibounded} allows us to prove that whenever
$\nu_k(V)<0$, then it is attained, thus providing a {\em bound state}
in $L^2(\R^N)$ with negative energy.

\begin{Corollary}\label{c:speck}
If $V\in{\mathcal V}$ and
$\nu_k(V)<0$, 
then $\nu_k(V)$ is attained.
\end{Corollary}

\begin{pfn}{Proposition \ref{p:spe}.2}
Since operators $-\D-V$,
$V\in{\mathcal V}$, are $L^{N/2}$-perturbations of positive operators
(see Lemma \ref{l:semibounded}),
from the Cwikel-Lieb-Rosenblum inequality (\cite{Cwikel, Lieb,
  Rosenblum}) it follows that the number of negative
eigenvalues is finite. Hence the conclusion follows from Corollaries
\ref{c:spec} and~\ref{c:speck}.  
\end{pfn}

\subsection{Eigenvalues at the bottom of the essential
spectrum}\label{sec:bottom-ess-spectr} 

\noindent We now mean to study the nature of the bottom of essential
spectrum of operators $L_{\lambda_1,\dots,\lambda_k,a_1,\dots,a_k}$ defined in
(\ref{eq:17}).  
More precisely, when  the values of  $\lambda_i$'s admit both 
configurations of poles corresponding to 
negative quadratic forms and configurations corresponding to positive
quadratic forms, we will provide a necessary and sufficient condition
on the masses of singularities 
for the existence of a configuration of $a_i$'s admitting a bound
state with null energy. 
 
Let $(\lambda_1,\dots,\lambda_k)\in\R^k$ fixed. We denote as $\Sigma$
the set of colliding configurations, 
namely 
$$
\Sigma:=\{(a_1,\dots,a_k)\in\R^{Nk}:\ a_i=a_j\
\text{for some }i\neq j\}. 
$$
For any
$\boldsymbol{a}=(a_1,\dots,a_k)\in\R^{Nk}\setminus\Sigma$,   
we introduce the following notation
$$
\mu_{\boldsymbol{a}}:=\inf_{u\in\Di\setminus\{0\}}\frac{
Q_{\lambda_1,\dots,\lambda_k,a_1,\dots,a_k}(u)}
{\|u\|^2_{\Di}}.
$$ 
The following result is a direct corollary of Lemma \ref{l:cneg}.

\begin{Corollary}\label{c:cont}
For $(\lambda_1,\dots,\lambda_k)\in (-\infty,(N-2)^2/4)^k$, let 
$\boldsymbol{a}_n\in\R^{Nk}$  be a
sequence of configurations converging to
$\boldsymbol{a}\in\R^{Nk}\setminus\Sigma$. Then
$\lim\limits_{n\to\infty}\mu_{\boldsymbol{a}_n}= \mu_{\boldsymbol{a}}$.
\end{Corollary}

\noindent Let us denote
$$
A^+\!\!:=\{\boldsymbol{a}\in\R^{Nk}\setminus \Sigma:\
\mu_{\boldsymbol{a}}>0\},\ 
A^-\!\!:=\{\boldsymbol{a}\in\R^{Nk}\setminus \Sigma:\
\mu_{\boldsymbol{a}}<0\}, \  \text{and}\  A^0\!\!:=\{\boldsymbol{a}\in\R^{Nk}\setminus \Sigma:\
\mu_{\boldsymbol{a}}=0\}.
$$
From Corollary \ref{c:cont}, it follows that $A^+$ and $A^-$ are open
sets. Hence, whenever both  $A^+$ and $A^-$ are nonempty, the set
$A^0$ is nonempty and disconnects $\R^{Nk}\setminus \Sigma$.

\bigskip\noindent

\begin{pfn}{Theorem \ref{t:bottom}}
  Let us assume that the $\lambda_i$'s satisfy (\ref{eq:19}) and
  (\ref{eq:72}).  From Theorem \ref{t:mainresult} and (\ref{eq:19}),
  there exists a configuration of poles
  $\boldsymbol{a^+}=(a_1^+,\dots,a_k^+)$ such that
  $\mu_{\boldsymbol{a^+}}>0$. On the other hand, in \cite[Proposition
  1.2]{FT} it is proved that, if $\sum_{i=1}^k\lambda_i^+>\frac{(N-2)^2}4$,
  then it is possible to find a configuration of poles
  ${\boldsymbol{a^-}}=(a_1^-,\dots,a_k^-)$ such that
  $\mu_{\boldsymbol{a^-}}<0$. It is worth noticing that, from the
  proofs of Theorem \ref{t:mainresult} and \cite[Proposition 1.2]{FT},
  there easily results that $\boldsymbol{a^+}$ and $\boldsymbol{a^-}$
  can be chosen to be collisionless, i.e.
  $\boldsymbol{a^+},\boldsymbol{a^-}\in\R^{Nk}\setminus\Sigma$.  From
  Corollary \ref{c:cont},
  ${\boldsymbol{a}}\mapsto\mu_{\boldsymbol{a}}>0$ is continuous on
  $\R^{Nk}\setminus\Sigma$, which is a connected open subset of
  $\R^{Nk}$. Therefore there exist
  ${\boldsymbol{a}}=(a_1,\dots,a_k)\in \R^{Nk}\setminus\Sigma$ such
  that $\mu_{\boldsymbol{a}}=0$.  From Proposition \ref{p:attai}, it
  follows that $\mu_{\boldsymbol{a}}=0$ is attained by some $u\in\Di$
  weakly solving in $\Di\setminus\{0\}$
\begin{equation}\label{eq:70}
-\D u(x)-\sum_{i=1}^k\frac{\lambda_i}{|x-a_i|^2}\,u(x)=0.
\end{equation}
By evenness we can assume $u\geq0$, while the
Strong Maximum Principle and standard regularity theory ensure that
$u$ is smooth and strictly positive in $\R^N\setminus\{a_1,\dots,a_k\}$.
Lemma \ref{l:2} yields a precise estimate of the decay of
$u$ at infinity, i.e. $u(x)\sim |x|^{-(N-2-a_{\lambda_{\infty}})}$ as
$|x|\to\infty$, where $\lambda_{\infty}=\sum_{i=1}^k\lambda_i$. As a consequence
$$
u\in L^2(\R^N)\quad\text{if and only if}\quad
\sum_{i=1}^k\lambda_i<\frac{(N-2)^2}4-1. 
$$
Hence, under assumption (\ref{eq:72}), any function $u$ attaining
$\mu_{\boldsymbol{a}}$ provides an eigenfunction of
the Schr\"o\-dinger operator $L_{\lambda_1,\dots,\lambda_k,a_1,\dots,a_k}$
associated to the null eigenvalue. 

Let us now prove the necessity of condition (\ref{eq:72}). If 
$\sum_{i=1}^k\lambda_i\geq\frac{(N-2)^2}4-1$,  Lemma \ref{l:2} implies
that, for any $\boldsymbol{a}\in \R^{Nk}$, 
(\ref{eq:70}) cannot have any nontrivial nonnegative solution in $H^1(\R^N)$.
On the other hand, if $\sum_{i=1}^k\lambda_i^+\leq\frac{(N-2)^2}4$, we
distinguish two cases:
\begin{description}
\item[Case 1] $\sum_{i=1}^k\lambda_i^+<\frac{(N-2)^2}4$. In this case, \cite[Proposition
1.2]{FT} ensures that $\mu_{\boldsymbol{a}}>0$ for any
$\boldsymbol{a}\in\R^{Nk}\setminus\Sigma$ and hence the only
$\Di$-solution to (\ref{eq:70}) is the null one.
\item[Case 2] $\sum_{i=1}^k\lambda_i^+=\frac{(N-2)^2}4$. In this case,
assumption (\ref{eq:19}) implies that there exists at least one
index $i$ such that $\lambda_i<0$. Arguing by contradiction, assume that
(\ref{eq:70})  admits a nontrivial $\Di$-solution $u$ to
(\ref{eq:70}) for some $(a_1,\dots,a_k)\in \R^{Nk}\setminus\Sigma$. We have that
\begin{align*}
0&=\bigg(1-\frac4{(N-2)^2}\bigg(\sum_{i=1}^k\lambda_i^+\bigg)\bigg)\int_{\R^N}|\n
u(x)|^2\,dx\\
&\leq \int_{\R^N}|\n
u(x)|^2\,dx-\sum_{i=1}^k\lambda_i^+\int_{\R^N}\frac{u^2(x)}{|x-a_i|^2}\,dx=
-\sum_{i=1}^k\lambda_i^-\int_{\R^N}\frac{u^2(x)}{|x-a_i|^2}\,dx<0,
\end{align*}
which is a contradiction.   
\end{description}  
In both cases, we have proved that, for any $(a_1,\dots,a_k)\in
\R^{Nk}\setminus\Sigma$,  (\ref{eq:70})  admits no nontrivial
$\Di$-solutions. In particular,  for any configuration of
singularities, $0$ is not an  eigenvalue of the Friedrichs extension of
$L_{\lambda_1,\dots,\lambda_k,a_1,\dots,a_k}$.
\end{pfn}

\appendix
\section*{Appendix}
\setcounter{section}{1}
\setcounter{equation}{0}
\setcounter{Theorem}{0}

We collect in this appendix some technical results used in the
paper. In the following lemma (which was needed in the proof of Theorem
\ref{t:self-adjo}), we extend to $L^2$ (not necessarily
bounded) functions a well-known property of
differentiability of {\em{Newtonian potentials}}, see \cite[Lemma 4.1,
p. 54]{GT}. 
\begin{Lemma}\label{l:green}
Let $\Omega\subset\R^N$ be a bounded smooth domain, $p\in \Omega$, $g\in
L^2(\Omega)$, $g$ smooth in $\Omega\setminus\{p\}$, and let
$u$ be the Newtonian potential of $g$, i.e.
$$
u(x)=\frac1{N(2-N)\omega_N}\int_{\Omega}\frac{g(y)}{|x-y|^{N-2}}\,dy,
\qquad x\in\R^N\setminus\{p\}.
$$
Then $u\in W^{1,q}(\R^N)$ for all
$q\in\big(\frac{N}{N-2},\frac{2N}{N-2}\big]$ and  the weak derivatives of
$u$ are given by 
$$
\frac{\partial u}{\partial
  x_i}(x)=\frac1{N\omega_N}
\int_{\Omega}\frac{g(y)(x_i-y_i)}{|x-y|^{N}}\,dy, 
\qquad x\in\R^N\setminus\{p\}. 
$$
\end{Lemma}

\begin{pf}
Let $\tilde g\in L^2(\R^N)$ be such that $\tilde g(y)=0$ in
$\R^N\setminus\overline\Omega$ and $\tilde g\big|_{\Omega}=g$. 
Note that $u=I_2(\tilde g)$, where $I_2(\tilde g)$ is the Riesz
potential defined by
$$
\big(I_2(\tilde g)\big)(x):=
\frac1{N(2-N)\omega_N}\int_{\R^N}\frac{\tilde
  g(y)}{|x-y|^{N-2}}\,dy.
$$
For any  $1< p<\min\{2,N/2 \}$, from \cite[Theorem 1,
p. 119]{stein} it is known that  $I_2$ is a linear bounded operator from
$L^{p}(\R^N)$ into $L^{( p N)/(N-2 p)}(\R^N)$.
It follows that $u\in L^q(\R^N)$ for all $q\in\big(\frac{N}{N-2},\infty\big)$
if $2<N\leq 4$, $q\in\big(\frac{N}{N-2},\frac{2N}{N-4}\big)$
if $N>4$.

Let $g_n\in  C^{\infty}_{\rm c}(\R^N)$ such that $\mathop{\rm
supp}g_n\subset\Omega$ and $g_n\to \tilde g$ in $L^p(\R^N)$ for all $p\in[1,2]$ and set
$u_n=I_2(g_n)$. Since for all  $1< p<\min\{2,N/2 \}$, $g_n\to \tilde g$  in $L^{ p}(\R^N)$, we
have that $u_n\to u$  in $L^{( p N)/(N-2 p)}(\R^N)$, i.e. $u_n\to u$
in $L^q(\R^N)$ for all $q\in\big(\frac{N}{N-2},\infty\big)$ 
if $2<N\leq 4$, $q\in\big(\frac{N}{N-2},\frac{2N}{N-4}\big)$
if $N>4$. By
\cite[Lemma 4.1, p. 54]{GT}, $u_n\in C^1(\R^N)$ and
$$
\frac{\partial u_n}{\partial
x_i}(x)=I_1^i(g_n):=\frac1{N\omega_N}
\int_{\R^N}\frac{g_n(y)(x_i-y_i)}{|x-y|^{N}}\,dy, 
\qquad x\in\R^N,\quad i=1,\dots,N. 
$$
From \cite[Theorem 1, p. 119]{stein}, $I_1^i$ are linear bounded operators from
$L^{p}(\R^N)$ into $L^{(p N)/(N-p)}(\R^N)$ for all $p\in(1,2]$. Hence,
for $i=1,\dots,N$, $\frac{\partial u_n}{\partial
  x_i}\to I_1^i(\tilde g)$ in $L^{(p N)/(N-p)}(\R^N)$ for all
$p\in(1,2]$, i.e.
$$
\frac{\partial u_n}{\partial
  x_i}\to 
\frac1{N\omega_N}
\int_{\R^N}\frac{\tilde g(y)(x_i-y_i)}{|x-y|^{N}}\,dy, 
\qquad \text{in }L^{q}(\R^N)\quad\text{for all }\frac{N}{N-1}<q\leq\frac{2N}{N-2}.
$$
Therefore for all $q\in\big(\frac{N}{N-2},\frac{2N}{N-2}\big]$, $u\in W^{1,q}(\R^N)$ and 
$$
\n u(x)=\frac1{N\omega_N}\int_{\Omega}\frac{g(y)(x-y)}{|x-y|^{N}}\,dy, 
\qquad x\in\R^N\setminus\{p\}. 
$$
The proof is thereby complete.
\end{pf}

\noindent The following lemma was used in the proof of Theorem
\ref{t:self-adjo}.

\begin{Lemma}\label{l:ode}
For $\bar s\in\R$, $\omega>0$, $b>0$, and $\a<0$, let $\psi\in C^1\big((-\infty,\bar s]\big)$ be the solution of the following Cauchy problem  
$$
\begin{cases}
\psi''(s)-\omega^2\,\psi(s)=b\, e^{2s}\,\psi(s),\\
\psi(\bar s)=0,\quad\psi'(\bar s)=\a.
\end{cases}
$$
Then 
$$
0\leq \psi(s)\leq -\frac{\a}{2\omega}\,e^{\omega\bar 
s}\exp\bigg(\frac{b}{4\omega}\,e^{2\bar s}\bigg)
\,e^{-\omega s}\quad\text{ for all }s\leq\bar s.
$$
\end{Lemma}

\begin{pf}
The initial conditions imply that $\psi$ is positive in a left
neighborhood of $\bar s$, whereas the equation forces the solution to be
convex wherever it is positive. As a consequence $\psi$ must be
strictly positive in $(-\infty,\bar s)$. We have that, for $s\leq\bar s$,
\begin{align*}
\psi(s)=e^{-\omega s}\bigg[-\frac{\a}{2\omega}e^{\omega\bar
s}-\frac{b}{2\omega}\int_{\bar s}^s e^{\omega t}\psi(t)e^{2t}\,dt\bigg]+
e^{\omega s}\bigg[\frac{\a}{2\omega}e^{-\omega\bar
s}+\frac{b}{2\omega}\int_{\bar s}^s e^{-\omega t}\psi(t)e^{2t}\,dt\bigg].
\end{align*} 
For any $\tau\geq0$, set $f(\tau):=e^{\omega(\bar s-\tau)}\psi(\bar s-\tau)$, hence
\begin{align*}
f(\tau)=-\frac{\a}{2\omega}e^{\omega\bar s}(1-e^{-2\omega\tau})+\frac{b}{2\omega}
\int_0^{\tau}\big(1-e^{2\omega(t-\tau)}\big)f(t)\,e^{2(\bar s-t)}\,dt.  
\end{align*}
Since the function $f\in C^1([0,+\infty))$ can be estimated as
$$
f(\tau)\leq -\frac{\a}{2\omega}\,e^{\omega\bar s}+\frac{b}{2\omega}\,
e^{2\bar s}\int_0^{\tau}f(t)\,e^{-2t}\,dt,
$$
the Gronwall's Lemma yields 
$$
f(\tau)\leq -\frac{\a}{2\omega}\,e^{\omega\bar s}\exp\bigg(\frac{b}{2\omega}\,
e^{2\bar s}\int_0^{\tau}e^{-2t}\,dt\bigg)\leq-\frac{\a}{2\omega}\,
e^{\omega \bar s}\exp\bigg(\frac{b}{4\omega}\,e^{2\bar s}\bigg)
$$
for all $\tau\geq0$, thus proving the required estimate.
\end{pf}

We now give a result of continuity of Hardy integrals with respect to
poles which was used in the proof of Lemma \ref{l:cneg}.  
\begin{Lemma}\label{l:cf}
For any $a\in\R^N$, $r>0$, and $u\in\Di$, there holds
$$
\lim_{y\to a}\int_{B(y,r)}\frac{u^2(x)}{|x-y|^2}\,dx=\int_{B(a,r)}\frac{u^2(x)}{|x-a|^2}\,dx.
$$
\end{Lemma}
\begin{pf}
For any $u\in \Di$, $u\geq 0$ a.e., we consider the Schwarz symmetrization
of $u$ defined~as
\begin{equation}\label{eq:67}
u^*(x):=\inf\big\{t>0:\ \big|\{y\in\R^N:\
u(y)>t\}\big|\leq \omega_N|x|^N\big\} 
\end{equation}
where $|\cdot|$ denotes the Lebesgue measure of $\R^N$ and
$\omega_N$ is the volume of the standard unit $N$-ball. For any $\Omega\subset\R^N$ 
and for any $u\in\Di$, let $\Omega^*=B\big(0,(|\Omega|/\omega_N)^{1/N}\big)$ and 
$|u|^*$ denote the Schwarz symmetrization of $|u|$, see \eqref{eq:67}. From
\cite[Theorem 21.8]{Willem} and since $\big(1/|x-y|\big)^*=1/|x|$, for any $y\in\R^N$,
it follows that
\begin{equation}\label{eq:68}
\int_{\Omega\cap B(y,r)}\frac{u^2}{|x-y|^2}\,dx\leq
\int_{\Omega^*\cap B(0,r)}\frac{(|u|^*)^2}{|x|^2}\,dx.
\end{equation}
Let $u\in\Di$. It is easy to see that 
$$
\frac{u^2\alchi_{B(y,r)}}{|x-y|^2}\quad\text{converges to}\quad
\frac{u^2\alchi_{B(a,r)}}{|x-a|^2} \quad\text{a.e. in }\R^N\text{ as }
y\to a. 
$$
Moreover, from (\ref{eq:68}), it follows that 
the family of functions $\Big\{\frac{u^2\alchi_{B(y,r)}}{|x-y|^2}:y\in\R^N\Big\}$ is equi-integrable. Hence Vitali's convergence Theorem allows to conclude.
\end{pf}


\end{document}